\documentclass[10pt]{amsart}
\usepackage{amsmath}
\usepackage{amscd,amsthm,amssymb,amsfonts}
\usepackage{mathrsfs}
\usepackage{dsfont}
\usepackage{stmaryrd}
\usepackage{euscript}
\usepackage{expdlist}
\usepackage{enumerate}

\input xy
\xyoption{all}
\usepackage[OT2,T1]{fontenc}

\theoremstyle{plain}
\newtheorem{thm}{Theorem}[section]

\newtheorem*{thm*}{Theorem}

\newtheorem{lm}[thm]{Lemma}
\newtheorem{cor}[thm]{Corollary}
\newtheorem*{cor*}{Corollary}

\newtheorem{prop}[thm]{Proposition}
\newtheorem*{conj*}{Conjecture}

\theoremstyle{remark}
\newtheorem*{remark}{Remark}
\newtheorem*{thank}{Acknowledgments}

\theoremstyle{definition}
\newtheorem*{defn*}{Definition}

\newtheorem{defn}[thm]{Definition}

\newcommand{\nc}{\newcommand}

\newcommand{\beq}{\begin{equation}}
\newcommand{\eeq}{\end{equation}}
\newcommand{\bpmx}{\begin{pmatrix}}
\newcommand{\epmx}{\end{pmatrix}}
\newcommand{\bbmx}{\begin{bmatrix}}
\newcommand{\ebmx}{\end{bmatrix}}
\newcommand{\wh}{\widehat}

\newcommand{\beqcd}[1]{\begin{equation*}\label{#1}\tag{#1}}
\newcommand{\eeqcd}{\end{equation*}}

\numberwithin{equation}{section}
\newenvironment{mylist}{
  \begin{enumerate}{}{%
      \setlength{\itemsep}{5pt} \setlength{\parsep}{0in}
      \setlength{\parskip}{0in} \setlength{\topsep}{0in}
      \setlength{\partopsep}{0in}
      \setlength{\leftmargin}{0.17in}}}{\end{enumerate}}

\def\makeop#1{\expandafter\def\csname#1\endcsname
  {\mathop{\rm #1}\nolimits}\ignorespaces}

\makeop{Hom}   \makeop{End}   \makeop{Aut}   
\makeop{Pic} \makeop{Gal}       \makeop{Div} \makeop{Lie}
\makeop{PGL}   \makeop{Corr} \makeop{PSL} \makeop{sgn} \makeop{Spf}
 \makeop{Tr} \makeop{Nr} \makeop{Fr} \makeop{disc}
\makeop{Proj} \makeop{supp} \makeop{ker}   \makeop{Im} \makeop{dom}
\makeop{coker} \makeop{Stab} \makeop{SO} \makeop{SL} \makeop{SL}
\makeop{Cl}    \makeop{cond} \makeop{Br} \makeop{inv} \makeop{rank}
\makeop{id}    \makeop{Fil} \makeop{Frac}  \makeop{GL} \makeop{SU}
\makeop{Trd}   \makeop{Sp} \makeop{Tr}    \makeop{Trd} \makeop{Res}
\makeop{ind} \makeop{depth} \makeop{Tr} \makeop{st} \makeop{Ad}
\makeop{Int} \makeop{tr}    \makeop{Sym} \makeop{can} \makeop{SO}
\makeop{torsion} \makeop{GSp} \makeop{Tor}\makeop{Ker} \makeop{rec}
\makeop{Ind} \makeop{Coker}
 \makeop{vol} \makeop{Ext} \makeop{gr} \makeop{ad}
 \makeop{Gr}\makeop{corank} \makeop{Ann}
\makeop{Hol} 
\makeop{Fitt} \makeop{Mp} \makeop{CAP}

\def\Isom{\ul{\mathrm{Isom}}}

\DeclareMathAlphabet{\mathpzc}{OT1}{pzc}{m}{it}
\DeclareSymbolFont{cyrletters}{OT2}{wncyr}{m}{n}
\DeclareMathSymbol{\SHA}{\mathalpha}{cyrletters}{"58}

\def\makebb#1{\expandafter\def
  \csname bb#1\endcsname{{\mathbb{#1}}}\ignorespaces}
\def\makebf#1{\expandafter\def\csname bf#1\endcsname{{\bf
      #1}}\ignorespaces}
\def\makegr#1{\expandafter\def
  \csname gr#1\endcsname{{\mathfrak{#1}}}\ignorespaces}
\def\makescr#1{\expandafter\def
  \csname scr#1\endcsname{{\EuScript{#1}}}\ignorespaces}
\def\makecal#1{\expandafter\def\csname cal#1\endcsname{{\mathcal
      #1}}\ignorespaces}

\def\doLetters#1{#1A #1B #1C #1D #1E #1F #1G #1H #1I #1J #1K #1L #1M
                 #1N #1O #1P #1Q #1R #1S #1T #1U #1V #1W #1X #1Y #1Z}
\def\doletters#1{#1a #1b #1c #1d #1e #1f #1g #1h #1i #1j #1k #1l #1m
                 #1n #1o #1p #1q #1r #1s #1t #1u #1v #1w #1x #1y #1z}
\doLetters\makebb   \doLetters\makecal  \doLetters\makebf
\doLetters\makescr
\doletters\makebf   \doLetters\makegr   \doletters\makegr

\normalsize

\makeop{Ram} \makeop{Rep} \makeop{mass}

\makeop{Bl}
\def\abs#1{\left|#1\right|}

\def\Qbar{\ol{\Q}}

\def\rmN{{\mathrm N}}

\def\cA{{\mathcal A}}  

\def\cE{{\mathcal E}}
\def\cF{{\mathcal F}}  
\def\cG{{\mathcal G}}
\def\cL{{\mathcal L}}

\def\cK{{\mathcal K}}  
\def\cM{\mathcal M}
\def\cR{{\mathcal R}}
\def\cO{\mathcal O}

\def\cf{{\mathcal f}}
\def\cW{{\mathcal W}}

\def\cT{\mathcal T}

\def\cU{\mathcal U}

\def\EucO{{\EuScript O}}

\def\bfc{\mathbf c}

\def\bff{\mathbf f}
\def\bdh{\mathbf h}

\def\bdc{\mathbf c}

\def\bftheta{\boldsymbol{\theta}}

\def\bbI{\mathbb I}

\newcommand{\Z}{\mathbf Z}
\newcommand{\Q}{\mathbf Q}
\newcommand{\R}{\mathbf R}
\newcommand{\C}{\mathbf C}
\newcommand{\A}{\mathbf A}    
\newcommand{\F}{\mathbb F}

\def\frakc{{\mathfrak c}}

\def\frakr{\mathfrak r}

\def\frakm{\mathfrak m}

\def\frakl{\mathfrak l}

\def\frakF{{\mathfrak F}}

\def\frakC{{\mathfrak C}}

\def\frakX{\mathfrak X}

\def\frakN{\mathfrak N}
\def\il{\mathfrak i\frakl}

\def\wbar{\bar{w}}

\def\wbar{\bar{w}}

\def\imply{\Rightarrow}

\def\ol{\overline}  \nc{\opp}{\mathrm{opp}} \nc{\ul}{\underline}

\newcommand{\MX}[4]{\begin{bmatrix}
{#1}& {#2}\\
{#3}&{#4}\end{bmatrix} }

 \newcommand{\DII}[2]{\begin{bmatrix}{#1}&0
 \\0&{#2}\end{bmatrix}}

\def\cf{\mbox{{\it cf.} }}

\def\can{{can}}

\def\uf{\varpi} 
\def\Abs{{|\!\cdot\!|}} 

\def\Sg{{\varSigma}}  
\def\ndivides{\nmid}
\def\ndivide{\nmid}
\def\x{{\times}}

\def\onehalf{{\frac{1}{2}}}

\def\iso{\simeq}
\def\con{\equiv}

\newcommand\stt[1]{\left\{#1\right\}}

\def\lam{\lambda}

\def\sg{\sigma}
\def\vp{\varphi}

\def\AFf{\A_{\cF,f}}

\def\AF{\A_\cF}

\def\setp{{(p)}}

\newcommand{\powerseries}[1]{\llbracket{#1}\rrbracket}

\renewcommand\pmod[1]{\,(\mbox{mod }{#1})}
\renewcommand\mod[1]{\,\mbox{mod }{#1}}

\usepackage[utf8x]{inputenc}
\usepackage{amsmath}
\usepackage{amssymb}
\usepackage{amsthm}
\usepackage{amsmath,amscd}
\usepackage[small,nohug,heads=vee]{diagrams}
\diagramstyle[labelstyle=\scriptstyle]
\title[On the non-triviality of the $p$-adic Abel--Jacobi image of generalised Heegner cycles]{On the non-triviality of the $p$-adic Abel--Jacobi image of generalised Heegner cycles modulo $p$,\\
 I: Modular curves}
\author{Ashay A. Burungale}
\date{\today}
\address{ 
Institute for Advanced Study, 1 Einstein Drive, Princeton, NJ 08540. \\
\newline
California Institute of Technolology,
1200 E California Blvd, Pasadena CA 91125. 
}
\email{ashayburungale@gmail.com}
\def\Csplit{\frakF}

\def\opcpt{K}
\def\lsgN{K^n_1}
\def\Section{\phi_{\chi,s,v}}

\def\cmptv{\varsigma_v}

\begin{document}
\maketitle

\begin{abstract} 
Generalised Heegner cycles are associated with a pair of an elliptic newform and a 
Hecke character over an imaginary quadratic extension $K/\Q$. 
Let $p$ be an odd prime split in $K/\Q$ and $\ell\neq p$ an odd unramified prime. 
We prove the non-triviality 
of the $p$-adic Abel--Jacobi image of generalised Heegner cycles 
modulo $p$ 
over the $\Z_\ell$-anticylotomic extension of $K$. 
The result is an evidence for the refined Bloch--Beilinson and the Bloch--Kato conjecture. 
In the case of weight two and $\ell$ an ordinary prime, it provides a non-trivial refinement of the results 
of Cornut and Vatsal on Mazur's conjecture regarding the non-triviality of Heegner points over the $\Z_\ell$-anticylotomic extension of $K$.
In the case of weight two and $\ell$ a supersingular prime, it settles Mazur's conjecture earlier known just for $\ell$ ordinary.\\
\end{abstract}
\tableofcontents

\def\ZZbox{\Z_{(\Box)\,}}
\def\cAbox{\cA^{(\Box)}_{K,\bdc}}
\def\cAboxn{\cA^{(\Box)}_{K,\bdc,n}}
\def\Zhatbox{\widehat{\Z}^{(\Box)}}
\def\qchKF{\tau_{\cK/\cF}}
\def\opcpt{K}
\def\sh{Sh}
\def\opn{K^n}
\def\lp{j}
\def\lpp{\eta^{(p)}}
\def\lsgN{\opcpt_1^n}
\def\Om{\boldsymbol{\omega}}
\def\wt{k}
\def\skewhf{\vartheta}
\def\Fv{F}
\def\Kv{E}
\def\OFp{\cO_{\cF_p}}
\def\ads{\chi}
\def\tbar{\ol{t}}
\def\wbar{\ol{w}}
\def\Section{\vp}
\def\adelef{\A_{\cF,f}}
\def\cmpt{\varsigma}
\def\cmptv{\varsigma_v}
\def\skewhf{\vartheta}
\def\OK{\cO_\cK}
\def\cmJ{\boldsymbol{J}}
\def\OFv{\cO_F}
\def\OKv{\cO_E}
\def\OF{O}
\def\LR{\cR}
\def\addchar{\psi}
\def\cK{E}
\def\cF{F}
\def\Op{O_p}
\def\torsbgp{\cU_p}
\def\loca{\mathrm{l.a.}}
\def\baseR{R}
\def\Prd{\boldsymbol{P}}
\def\wtsp{\frakX^{\mathrm{crit}}_{p}}
\def\IwGamma{\Gamma^-}
\def\AF{\A}
\def\AFf{\AF_f}
\def\CMP{\vartheta}
\def\OKbasis{\bftheta}
\newcommand\localW[1]{W_{#1}}
\newcommand\localK[1]{\bfa_{#1}}
\newcommand\localP[2]{\Prd(#1,#2)}
\def\rmN{\mathrm N}
\def\Sgbar{\ol{\Sg}}
\renewcommand\MX[4]{\begin{pmatrix}#1&#2\\
#3&#4\end{pmatrix}}
\renewcommand\DII[2]{\begin{pmatrix}#1&\\
&#2\end{pmatrix}}
\def\Cinert{\frakN^-}
\def\Csplit{\frakF}
\def\Cram{D_{\cK/\cF}}
\def\wt{k}
\def\plideal{\frakc}
\def\dt{dt}
\def\Bad{\frakr}
\def\CMring{{\cW_p}}
\def\#{\sharp}

\section{Introduction} 
\noindent
One expects the existence of a commutative $p$-adic L-function $L_p(\cM)$ associated to a $p$-ordinary critical motive $\cM$ 
over a tower of number fields with the Galois group being a commutative $p$-adic Lie group $\Gamma$, characterised by
an interpolation of the $p$-stabilised 
critical L-values of $\cM$ twisted by a dense subset of characters of $\Gamma$. Certain values of the $p$-adic L-function outside the range of interpolation 
are expected to encode deep $p$-adic arithmetic information regarding $\cM$. 
One can ask for the non-triviality of the value of $L_p(\cM \otimes \nu)$ at a fixed character outside the range of interpolation, 
as $\nu$ varies over a family of $\ell$-power conductor characters for a prime $\ell \neq p$. 
If one has a strategy of approaching the non-triviality of 
critical L-values of $\cM$ twisted by $\nu$ as $\nu$ varies, it might apply for the former non-triviality as well.\\
\\
Based on the classical Waldspurger formula (\cf \cite{W}), an anticyclotomic $p$-adic L-function associated with a Rankin--Selberg convolution of an elliptic modular 
form and a theta series is constructed in \cite{BDP1}, \cite{Br} and \cite{Hs3}. 
In \cite{BDP1}, values of the $p$-adic L-function at certain characters outside the range of interpolation are 
shown to be related to the $p$-adic Abel--Jacobi image of generalised Heegner cycles associated to the convolution. 
This can be considered as a $p$-adic Waldspurger formula (\cf \cite{LZZ}). 
The construction of the cycle is due to Bertolini--Darmon--Prasanna and generalises the one of classical Heegner cycles due to Schoen and Nekov\'{a}\v{r}. 
The cycle lives in a middle dimensional Chow group of a fiber product of a Kuga--Sato variety arising from a modular curve 
and a self-product of a CM elliptic curve. 
In the case of weight two, the cycle coincides with a Heegner point. 
Based on Hida's strategy, non-triviality of Rankin--Selberg L-values of the $l$-power conductor anticyclotomic character twists modulo $p$ 
 is shown in \cite{Br1} and \cite{Hs3}. 
In this article, we use the strategy to study non-triviality of the $p$-adic Abel-Jacobi image of generalised Heegner cycles modulo $p$ 
over the $\Z_l$-anticyclotomic extension. 
The non-triviality can be seen as an evidence for the Bloch--Beilinson and the Bloch--Kato conjecture as follows. 
In the setup, the Rankin--Selberg convolution is self-dual with root number $-1$. 
Accordingly, the Bloch--Beilinson conjecture predicts the existence of a non-torsion null-homologous cycle in the Chow realisation. 
In the setup, a natural candidate for a non-trivial null-homologous cycle is the generalised Heegner cycle. 
For a prime $p$,  the Bloch--Kato conjecture implies the non-triviality of the $p$-adic étale Abel--Jacobi image of the cycle. 
A natural question is to further investigate the non-triviality of the $p$-adic Abel--Jacobi image of the cycle.\\
\\
In the introduction, for simplicity we mostly restrict to the case of Heegner points.\\
\\
Let $\Qbar$ be an algebraic closure of $\Q$.
We fix two embeddings $\iota_{\infty}\colon \Qbar \to \C$ and $\iota_{p}\colon \Qbar \to \C_p$ 
for a prime $p$.
Let $v_p$ be the $p$-adic valuation induced by $\iota_p$ so that $v_p(p)=1$. 
Let $\mathfrak{m}_p$ be the maximal ideal of $\overline{\Z}_p$.
Unless otherwise stated, we suppose $p>3$.\\
\\
Let $K/\Q$ be an imaginary quadratic extension and $\cO$ the ring of integers. 
Let $-d_{K}$ be the discriminant. 
As $K$ is a subfield of the complex numbers, we regard it as a subfield of the algebraic closure $\Qbar$ via the embedding $\iota_{\infty}$. 
Let $G_{K}$ be the absolute Galois group $\Gal(\overline{\Q}/K)$. 
Let $c$ be the complex conjugation on $\C$ which induces the unique non-trivial element of $\Gal(K/\Q)$ via $\iota_{\infty}$. 
We assume the following:\\
\\
{(ord)}  \text{$p$ splits in $K$.}\\
\\
Let $\mathfrak{p}$ be a prime above $p$ in $K$ induced by the $p$-adic embedding $\iota_p$. For an integral ideal $\mathfrak{n}$ of $K$, 
we fix a decomposition $\mathfrak{n=n^{+}n^{-}}$ where $\mathfrak{n}^{+}$ (resp. $\mathfrak{n}^{-}$) is only divisible by split (resp. ramified or inert) primes in $K/\Q$. 
For a positive integer $m$, let $H_{m}$ be the ring class field of $K$ with conductor $m$ and 
$\cO_{m} = \Z+m\cO$ the corresponding order.
Let $H$ be the Hilbert class field.\\
\\
Let $N$ be a positive integer such that $p\ndivide N$. 
For $k\geq 2$, let $S_{k}(\Gamma_{0}(N),\epsilon)$ be the space of elliptic cusp forms of weight $k$, level $\Gamma_{0}(N)$ and neben-character $\epsilon$.
Let $f\in S_{2}(\Gamma_{0}(N),\epsilon)$ be a normalised newform. 
In particular, it is a Hecke eigenform with respect to all Hecke operators. 
Let $N_{\epsilon}|N$ be the conductor of $\epsilon$. 
Let $E_f$ be the Hecke field of $f$ and $\cO_{E_f}$ the ring of integers. 
Let $\Z_f$ be the order generated by the Hecke eigenvalues.
Let $\mathfrak{P}$ be a prime above $p$ in $E_f$ induced by the $p$-adic embedding $\iota_p$. 
Let $\rho_{f}:\Gal(\overline{\Q}/\Q)\rightarrow \GL_{2}(\cO_{E_{f,\mathfrak{P}}})$ be the corresponding $p$-adic Galois representation.\\
\\
We assume the following generalised Heegner hypothesis:\\
\\
{(Hg)} The integer ring $\cO$ contains an ideal $\mathfrak{N}$ of norm $N$ such that there exists an isomorphism 
$$\cO/\mathfrak{N} \iso \Z/N\Z .$$
From now, we fix such an ideal $\mathfrak{N}$. Let $\mathfrak{N}_{\epsilon}|\mathfrak{N}$ be the unique ideal of norm $N_{\epsilon}$.
In view of (Hg), note that $\frak{N}^{-}$ is only divisible by ramified primes in $K$.\\
\\
Let ${\bf{N}}$ denote the norm Hecke character over $\Q$ and ${\bf{N}}_{K}:={\bf{N}}\circ N_{\Q}^{K}$ the norm Hecke character over $K$.
For a Hecke character $\lam$ over $K$, let $\mathfrak{f}_\lam$ (resp. $\epsilon_\lam$) denote 
its conductor (resp. the restriction $\lam|_{{\bf{A}}_{\Q}^\times}$, where ${\bf{A}}_{\Q}$ denotes the adele ring over $\Q$).
We say that $\lam$ is central critical for $f$ if it is of infinity type $(j_{1},j_{2})$ with $j_{1}+j_{2}=2$ and $\epsilon_{\lam}=\epsilon {\bf{N}^{2}}$.\\
\\
Let $b$ be a positive integer prime to $pN$. 
Let $\Sg_{cc}(b,\mathfrak{N},\epsilon)$ be the set of Hecke characters $\lam$ such that:\\
\\
(C1) $\lam$ is central critical for $f$,\\
(C2) $\mathfrak{f}_{\lam}=b \cdot \mathfrak{N}_\epsilon$ and\\
(C3) we have $\epsilon_{q}(f,\lam^{-1})=1$ for $\epsilon_{q}(f,\lam^{-1})$ the local root number corresponding to Rankin--Selberg convolution of the pair $(f,\lam^{-1})$ for all finite primes $q$.\\
\\
Here we use the subscript `$cc$' as the complex L-function $L(s, f \times \lam)$ corresponding to the Rankin--Selberg convolution arising from the pair $(f,\lam)$ 
has $s=1$ as a central-critical point.\\
\\
Let $\chi$ be a finite order Hecke character such that $\chi{\bf{ N}}_{K} \in \Sg_{cc}(b,\mathfrak{N},\epsilon)$. In view of hypothesis (C3), we have $$\epsilon(f, \chi^{-1})=-1$$
for the global root number $\epsilon(f, \chi^{-1})$ of the Rankin--Selberg convolution.
Let $E_{f,\chi}$ be the finite extension of $E_f$ obtained by adjoining the values of $\chi$.\\
\\
Let $X_1(N)$ be the modular curve of level $\Gamma_1(N)$, $\infty$ the standard cusp $i\infty$ of $X_1(N)$ and $J_1(N)$ the corresponding Jacobian. 
Let $B_f$ be an abelian variety associated to $f$ by the Eichler--Shimura correspondence 
and $\Phi_{f}:J_1(N) \rightarrow B_f$ the associated surjective morphism. 
Up to isogeny, the uniqueness of the abelian variety follows from the assumption that $f$ is a newform. 
By possibly replacing $B_f$ with an isogenous abelian variety, we suppose
that $B_{f}$ has endomorphisms by the order $\Z_{f}$. 
Let $\omega_{f}$ be the differential form on $X_{1}(N)$ corresponding to 
$f$. We use the same notation for the corresponding $1$-form on $J_1(N)$. 
Let $\omega_{B_{f}} \in \Omega^{1}(B_{f}/E_{f})^{\Z_{f}}$ be the unique $1$-form such that $\Phi_{f}^{*}(\omega_{B_{f}})=\omega_{f}$. 
Here $\Omega^{1}(B_{f}/E_{f})^{\Z_{f}}$ denotes the subspace of $1$-forms given by 
$$
\Omega^{1}(B_{f}/E_{f})^{\Z_{f}}=\bigg{\{}\omega \in \Omega^{1}(B_{f}/E_{f})|[\lam]^{*} \omega = \lam \omega, \forall \lam \in \Z_{f}\bigg{\}}.
$$
\noindent\\
Recall that $b$ is a positive integer prime to $N$. Let $A_{b}$ be an elliptic curve with endomorphism ring $\cO_{b}=\Z+b\cO$, defined over the ring class field $H_b$. 
Let $t$ be a generator of $A_{b}[\mathfrak{N}]$. We thus obtain a point $x_{b}=(A_{b},A_{b}[\mathfrak{N}],t) \in X_{1}(N)(H_{bN})$. 
Let $\Delta_{b}=[A_{b},A_{b}[\mathfrak{N}],t]- (\infty) \in J_{1}(N)(H_{bN})$ be the corresponding Heegner point on the modular Jacobian. 
We regard $\chi$ as a character $\chi: \Gal(H_{bN}/K) \rightarrow E_{f,\chi}$. 
Let $H_{\chi}$ be the abelian extension of $K$ cut out by the character $\chi$.
To the pair $(f,\chi)$, we associate the Heegner point $P_{f}(\chi)$ given by 
\beq
P_{f}(\chi)= \sum_{\sg \in \Gal(H_{bN}/K)} \chi^{-1}(\sg)\Phi_{f}(\Delta_{b}^{\sg}) \space \in B_{f}(H_{\chi}) \otimes_{\Z_{f}} E_{f,\chi}.
\eeq
To consider the non-triviality of the Heegner points $P_{f}(\chi)$ as $\chi$ varies, we can consider the non-triviality of the 
corresponding $p$-adic formal group logarithm. 
The restriction of the $p$-adic logarithm gives a homomorphism 
$$\log_{\omega_{B_{f}}}:  B_{f}(H_{\chi}) \rightarrow \C_{p}.$$ 
We extend it to $B_{f}(H_{\chi}) \otimes_{\Z_{f}} E_{f,\chi}$ by $E_{f,\chi}$-linearity.\\
\\
We now fix a finite order Hecke character $\eta$ such that $\eta {\bf{N}}_{K} \in \Sg_{cc}(c,\mathfrak{N},\epsilon)$, for some $c$. 
Such an $\eta$ indeed exists for a sufficiently large integer $c$ (\cite[\S2.5]{BD}).
For $v|c^{-}$, let $\Delta_{\eta,v}$ be the finite group $\eta(\cO_{K_{v}}^{\times})$. 
Here $\cO_{K_{v}}$ denotes the integer ring of the local field $K_{v}$. 
Let $\ell\neq p$ be an odd prime unramified in $K$ and prime to $cN$. 
Let $H_{cN\ell^{\infty}} = \bigcup_{n\geq 0} H_{cN\ell^{n}}$. 
Let $G_{n}=\Gal(H_{cN\ell^{n}}/K)$
and $\Gamma_{\ell}=\varprojlim G_{n}$. 
Let $\mathfrak{X}_{\ell}$ denote the set of $\ell$-power order anticyclotomic characters of $\Gamma_\ell$. 
As $\nu \in \mathfrak{X}_{\ell}$ varies, the global root number $\epsilon(f,\eta^{-1}\nu^{-1})$ still equals $-1$. 
We consider the earlier setup of Heegner points with $\lambda$ being $\eta\nu$. 
A basic question may then be the mod $p$ non-triviality of Heegner points $P_{f}(\eta\nu)$. It is in turn closely related to $p$-indivisibility of corresponding normalised $p$-adic logarithms.\\
\\
Our result is the following.\\
\\
\\
{{\bf{Theorem A}}}. 
Let $p$ be an odd prime and $N$ a positive integer such that $p \ndivide N$. 
Let $K$ be an imaginary quadratic field satisfying (ord) and (Hg) with $\mathfrak{p}$ a prime above $p$ determined via an initial embedding $\iota_p: \overline{\Q} \hookrightarrow \C_p$. 
Let $f\in S_{2}(\Gamma_{0}(N),\epsilon)$ be a normalised newform with $T_{p}$-eigenvalue $\bfa_{p}(f)$ and $\eta$ a finite order Hecke character such that $\eta {\bf{N}}_{K} \in \Sg_{cc}(c,\mathfrak{N},\epsilon)$ 
for an integer $c$ prime to $pN$.
Suppose the following holds.\\
\\
(1). The residual representation $\rho_{f}|_{G_{K}} \mod{\mathfrak{m}_{p}}$ is absolutely irreducible,\\
(2). $N^-$ is square-free and\\
(3). $p \ndivide \prod_{v|c^{-}}|\Delta_{\eta,v}|$ for $\Delta_{\eta,v}$ the finite group $\eta(\cO_{K_{v}}^{\times})$.\\
\\
Let $\ell\neq p$ be an odd prime unramified in $K$ and prime to $cN$. 
Let $\mathfrak{X}_\ell$ be the set of $\ell$-power order anticyclotomic Hecke characters over $K$ as above. 
Then, for all but finitely many $\nu \in \mathfrak{X}_{\ell}$ we have
$$v_{p}\big{(} \cE(f,\eta\nu)\log_{\omega_{B_{f}}}(P_{f}(\eta \nu))\big{)}=0.$$
Here $\cE(f,\eta\nu)=1 - (\eta\nu)^{-1}(\overline{\frak{p}})\bfa_{p}(f)p^{-1}+(\eta\nu)^{-2}(\overline{\frak{p}})\epsilon(p)p^{-1}$.
\\
\\
\\
We have the following immediate corollary.\\
\\
\\
{{\bf{Corollary A}}}. 
Let $f\in S_{2}(\Gamma_{0}(N),\epsilon)$ be a normalised newform and $\eta$ a finite order Hecke character such that $\eta {\bf{N}}_{K} \in \Sg_{cc}(c,\mathfrak{N},\epsilon)$ for an integer $c$ prime to $pN$. 
In addition to (Hg), suppose that\\
\\
(1). $f$ does not have CM over $K$ and\\
(2). $N^-$ is square-free.\\
\\
Then,  
for all but finitely many $\nu \in \mathfrak{X}_{\ell}$ the Heegner points $P_{f}(\eta\nu)$ are non-zero in $B_{f}(H_{\eta\nu}) \otimes_{\Z_{f}} E_{f,\eta\nu}$.\\
\\
\\
\begin{proof}
It suffices to show the existence of a prime $p$ satisfying the hypotheses of Theorem A. 
In view of (1), the existence follows.\\
\end{proof}
\noindent\\
For analogous non-triviality of the $p$-adic Abel--Jacobi image of generalised Heegner cycles modulo $p$, we refer to \S4.3 (\cf Theorem 4.6).\\
\\
We now describe strategy of the proof. Some of the notation used here is not followed in the rest of the article. 
Our approach is based on a strategy of Hida. 
To slightly strengthen the strategy, we first define anticyclotomic toric periods $P_{g,\lam}(\nu,n)$ associated to a pair $(g,\lam)$. 
Here $g$ is a $p$-adic elliptic modular form, $\lam$ a Hecke character and $\nu$ a Hecke character as above factoring through $G_n$. 
We then state a non-triviality theorem providing a criterion for non-triviality of the toric periods modulo $p$ as long as 
$g$ is a $p$-integral nearly holomorphic elliptic modular form defined over a number field. 
The criterion involves non-triviality of a certain toric form associated to the pair modulo $p$. 
The theorem follows from Hida's proof of the non-triviality of an anticyclotomic modular measure modulo $p$ in \cite{Hip}. 
The latter is based on Chai's theory of Hecke stable subvarieties of a mod $p$ Shimura variety (\cf \cite{Ch3}).
Let $d$ be the Katz $p$-adic differential operator and $f^{(p)}$ the $p$-depletion of $f$. 
For the pair $(d^{-1}(f^{(p)}), \eta)$, 
the toric period $P_{d^{-1}(f^{(p)}), \eta}(\nu,n)$ essentially equals $\log_{\omega_{B_{f}}}(P_{f}(\eta \nu))/p$. 
This is essentially the $p$-adic Waldspurger formula due to Bertolini--Darmon--Prasanna. 
As $d^{-1}(f^{(p)})$ is a weight zero $p$-adic elliptic modular form, the non-triviality theorem does not directly apply. 
However, we show that there exists a pair $(h,\chi)$ such that its toric periods are congruent modulo $p$ to the toric periods of interest and 
the non-triviality theorem also applies.
The nearly holomorphic elliptic modular form $h$ is an iterated image under the $p$-adic differential operator of $f$ and 
the congruence is a rather simple consequence of the $q$-expansion principle. 
Based on computation of Fourier coefficients in \cite{Hs3} for toric modular forms, 
we verify the criterion in the non-triviality theorem for the pair $(h,\chi)$.
The hypotheses (1) and (2) guarantee non-triviality of an underlying toric form modulo $p$.\\
\\
We would like to emphasise that the approach to non-triviality of arithmetic invariants modulo $p$ due to Hida (\cite{Hip}) and an automorphic computation due to Hsieh (\cite{Hs3}) play an essential role in the strategy after commencing with $p$-adic Waldspurger formula due to Bertolini--Darmon--Prasanna (\cite{BDP1}).\\
\\
We would like to remark that neither L-values nor $p$-adic L-functions appear explicitly in the statement or the proof of Theorem A. The above strategy can be also used to prove the non-triviality of 
the $p$-adic formal group logarithm arising from toric forms in \cite{Hs3} of the Heegner points on the modular Jacobian modulo $p$ .\\
\\
Without the hypotheses (1) and (2), Corollary A was conjectured by Mazur (\cf \cite{M}). 
Under classical Heegner hypotheses and $\ell$ an ordinary prime, the conjecture was independently proven by Cornut and Vatsal for (\cf \cite{C}, \cite{V}, \cite{V1} and also \cite{CV}). 
For a related work of Aflalo--Nekov\'{a}\v{r}, we refer to \cite{AN}.
As far as we know, Corollary A is a first result regarding Mazur's conjecture in the case $\ell$ a supersingular prime. The earlier results also did not allow the generalised Heegner hypotheses (Hg).
Our approach fundamentally differs from the one in \cite{C} and \cite{V}.
In particular, it uniformly treats the ordinary and supersingular case.  
It seems suggestive to compare our approach with the earlier ones. 
Here we only mention the following and invite the reader to further analogies (or refer to \cite{Bu1}).
In \cite{V1}, Jochnowitz congruence is a starting point. It reduces the non-triviality of the Heegner points to the non-triviality of the Gross points on 
a suitable definite Shimura ``variety". This in turn is studied via Ratner's theorem. In this article, the non-triviality is based on the modular curve itself. 
As indicated above, it relies on Chai's theory of Hecke stable subvarieties of a mod $p$ Shimura variety. 
In particular, the result of Cornut and Vatsal now fits in a framework due to Hida (\cf \cite{Hip}). 
An expectation along these lines was expressed in \cite{V2}.
Before the $p$-adic Waldspurger formula, 
the result and the framework appeared to be complementary. 
The formula also allows a rather smooth transition to the higher weight case as well.\\
\\ 
As far as we know, Theorem A is a first general result result regarding the non-triviality of the $p$-adic formal group logarithm of Heegner points modulo $p$. 
A related modulo $p$ non-triviality is also considered in \cite{V1}. 
As \cite{V1} treats the case when $p$ is inert in $K$, our result is complementary.\\
\\
As far as we know,  
higher weight analogue of Theorem A is a first general result regarding the non-triviality of the $p$-adic Abel--Jacobi 
image of generalised Heegner cycles. 
The only known result seems to be the non-triviality of several examples of such families in \cite{BDP3}. 
Howard proved a result towards an analogue of non-triviality of the étale Abel--Jacobi image of classical Heegner cycles in \cite{Ho}.\\ 
\\
Theorem A and its higher weight analogue have various arithmetic applications 
besides the Bloch--Beilinson and the Bloch--Kato conjecture. 
Consequences for the $p$-indivisibility of Heegner points when $f$ is $p$-ordinary and 
the refined Bloch--Beilinson conjecture are described in \S5 and the remark (1) following Theorem 4.6, respectively. 
The latter builds on a consideration in \cite{BDP3}. 
The non-triviality also implies that the Griffiths group of the fiber product of a Kuga--Sato variety arising from the Modular curve
and a self-product of a CM elliptic curve has infinite rank over the anticyclotomic extension. 
Arithmetic consequences of Corollary A are well documented in the literature, for example \cite{C} and \cite{V1}.\\
\\
In \cite{Bu}, we prove an analogue of the results for Shimura curves over the rationals in the case $\ell=p$. 
Recently, the $p$-adic Waldspurger formula has been generalised to modular forms on Shimura curves over a totally real field (\cf \cite{LZZ}).\\
\\
The article is perhaps a follow up to the articles \cite{BDP1} and \cite{Hip}. We refer to them for a general introduction. For a survey, we refer to \cite{Bu1}. For non-triviality based on Zariski density of CM points in characteristic $p$ or zero, we refer to \cite{BuHi1}, \cite{BuHi2} and \cite{BuTi}.\\
\\
The article is organised as follows. 
In \S2, we recall certain generalities about modular curves and modular forms.
In \S3, we describe underlying results regarding mod $p$ non-triviality of the anticyclotomic toric periods. 
In \S3.1, we firstly describe the notion of anticyclotomic toric periods associated to a class of modular forms 
and then discuss their non-triviality based on a linear independence due to Hida. 
In \S3.2, we describe the linear independence of mod $p$ modular forms arising from Chai's theory of Hecke-stable subvarieties of a mod $p$ Shimura variety. 
In \S4, we prove non-triviality of generalised Heegner cycles modulo $p$. 
For simplicity, we treat the case of Heegner points individually. 
In \S4.1, we first state the $p$-adic Waldspurger formula 
and then prove Theorem A. 
In \S4.2, we describe generalities regarding generalised Heegner cycles and state conjectures regarding their non-triviality.
In \S4.3, we prove the non-triviality of the $p$-adic Abel--Jacobi image of generalised Heegner cycles modulo $p$. 
In \S5, we deduce $p$-indivisibility of Heegner points based on the main result in \S4.1.
In the appendix, we describe generalities regarding toric modular forms which is a key automorphic ingredient in the non-vanishing.\\
\begin{thank} 
We are grateful to our advisor Haruzo Hida for persistent guidance and encouragement. 
During the preparation of the first draft, the author was a
graduate student in UCLA. We are grateful to friends back then for the distinctive
atmosphere. 

We thank Ming-Lun Hsieh for numerous instructive comments on the previous versions of the article and also for his assistance. 
We thank Francesc Castella for helpful conversations, particularly regarding the $p$-adic Waldspurger formula. 
We thank Barry Mazur, Dinakar Ramakrishnan, Christopher Skinner and Burt Totaro for instructive comments and encouragement. 
We also thank Miljan Brakocevic, Brian Conrad, Henri Darmon, Samit Dasgupta, Chandrashekhar Khare, Jan Nekov\'{a}\v{r}, Ye Tian, Vinayak Vatsal, 
Kevin Ventullo and Xinyi Yuan for helpful conversations about the topic. 

Finally, we are indebted to the referee. 
The current form of the article owes much to the detailed comments and suggestions of the referee. 
\end{thank}
\noindent\\
\\
\noindent {\bf{Notation}}
We use the following notation unless otherwise stated.\\
\\
For a number field $L$, let $\cO_L$ be the ring of integers, ${\bf{A}}_{L}$ the adele ring, ${\bf{A}}_{L,f}$ the finite adeles and ${\bf{A}}_{L,f}^{(\Box)}$ the finite adeles away from a finite set of places 
$\Box$ of $L$.  
Let $h_L$ be the set of finite places of $L$ and $h=h_\Q$. 
Let ${\bf{A}}$ denote the adeles over $\Q$.
For a fractional ideal $\mathfrak{a}$, let $\widehat{\mathfrak{a}}=\mathfrak{a}\otimes_{\Z} \widehat{\Z}$. 
For $a \in L$, let $\mathfrak{il}_{L}(a)=a(\cO_L \otimes_{\Z} \widehat{\Z}) \cap L$. 

Let $G_L$ be the absolute Galois group of $L$ and $G_{L}^{ab}$ the maximal abelian quotient. 
Let $\rec_{L}: {\bf{A}}_{L}^\times \rightarrow G_{L}^{ab}$ be the geometrically normalized reciprocity law.\\ 
\\
We refer to an elliptic modular form as a modular form. For a modular form $g$, let its Fourier expansion $g(q)$ at a cusp ${\bf{c}}$ be given by
$$ g(q) = \sum_{n \geq 0} \bfa_{n}(g,{\bf{c}})q^{n}.$$
\noindent For a positive integer $n$, let $\mu_{n}$ denote the set of primitive $n^{th}$-roots of unity in an algebraically closed field (as relevant).\\
\\
For $S \subset T$, let $\bbI_{S}: T \rightarrow \{0,1\}$ denote the indicator function of $S$.
\noindent\\
\\
\section{Modular curves and modular forms}
\noindent In this section, we recall certain generalities regarding modular curves and modular forms.\\
\\
\subsection{Modular curves}
\noindent In this subsection, we recall certain generalities regarding modular curves.\\
\\
In regards to the article, the section is preliminary. 
It briefly recalls the geometric theory of modular curves, $p$-ordinary Igusa tower and CM points on them.\\
\\
\subsubsection{Setup} In this subsection, we recall a basic setup regarding modular curves. 
We refer to \cite[\S2.2]{Hip} and \cite[\S2.1]{Hs2} for details.\\
\\
Let $G= \GL_{2/\Z}$ and 
$X (\iso (\C - \R))$ the corresponding Hermitian symmetric domain. 
The pair gives rise to a tower $(Sh_{K}=Sh_{K}(G,X))_K$ of quasi-projective smooth curves over $\Q$ indexed by open compact subgroups $K$ of $G(\A_{\Q,{f}})$. 
The pro-algebraic variety $Sh_{/\Q}$ is the projective limit of these curves. 
The complex points of the varieties are given as follows\\
\beq 
Sh_{K}(\C) =  G(\Q) \backslash X \times G(\A_{\Q,f}) / K , Sh(\C)= G(\Q) \backslash X \times G(\A_{\Q,f}) / \overline{Z(\Q)}.\\
\eeq
\noindent\\
Here $\overline{Z(\Q)}$ is the closure of the center $Z(\Q)$ in $G(\A_{\Q,f})$ under the adelic topology.\\
\\
Let us introduce some notation. Consider $V = \Q^2$ as a two dimensional vector space over $\Q$. Let $e_{1} = (1,0)$ and $e_{2}=(0,1)$. 
Let $\langle . , . \rangle: V \times V \to \Q$ be the $\Q$-bilinear pairing defined by $\langle e_{1},e_{2}\rangle =1$. 
Let $ \cL = \Z e_{1} \oplus_{\Z} \Z e_{2}$ be the standard lattice in V. 
For $g \in G(\Q)$, $g':= det(g) g^{-1}$ and 
consider the left action $gx := xg'$ with $x \in V$.\\
\\
Let $N$ be a positive integer, $h$ the set of finite places of $\Q$ as above and $v \in h$. 
Let 
$$U(N)=\bigg{\{} g \in G(\A_{\Q,f}) \big{|} g \equiv 1 \mod{N\cL}\bigg{\}},
K_{v}^{0} =\bigg{\{} g \in GL_{2}(\Q_{v}) \big{|} g( \cL \otimes_{\Z} \Z_{v}) = \cL \otimes_{\Z} \Z_{v} \bigg{\}}.
$$
\noindent\\
Let $p$ be an odd prime.
From now on, we consider only the decomposable open compact subgroups $K=\prod_{v} K_{v}$ of $G(\A_{\Q,f})$ for which $K_p$ equals $K_{p}^{0}$ and $U(N) \subset K$, for some prime-to-$p$ integer $N$. 
For a prime $\ell$ and a non-negative integer $n$, let
$$K_{0}(\ell)=\bigg{\{} g \in K \big{|} e_{2}g \in \Z e_{2} \mod{\ell \cL}\bigg{\}}, 
K_{1}^{n}=\bigg{\{} g \in K \big{|} g \equiv \MX{1}{*}{0}{1}\mod{p^n}\bigg{\}}.
$$
\noindent\\
\\
\subsubsection{$p$-integral model} 
In this subsection, we briefly recall generalities about a canonical $p$-integral smooth model of 
a modular curve. 
The description underlies the $p$-integral structure on the space of elliptic modular forms.
We refer to \cite[\S4.2]{Hi1} for details.\\
\\
Let the notation and hypotheses be as in the introduction and \S2.1. 
Modular Shimura variety $Sh_{/\Q}$ represents a functor classifying elliptic curves 
along with additional structure
(\cf\cite[\S4.2]{Hi1} and \cite{Sh1}). 
A $p$-integral 
interpretation of the functor leads to a $p$-integral smooth model of $Sh/G(\Z_p)_{/\Q}$.
The $p$-integral interpretation is given by\\
$$\cG^{(p)}: SCH_{/\Z_{(p)}} \to SETS$$
\beq S \mapsto \big{\{}(A,\eta^{(p)})_{/S}\big{\}}/ \sim .\eeq
Here,\\
\\
(PM1) $A$ is an elliptic curve over $S$.\\
(PM2) Let $\cT^{(p)}(A)$ be the prime-to-$p$ Tate module $\varprojlim_{p \ndivide N} A[N]$. 
$\eta^{(p)}$ is a prime-to-$p$ level structure given by a  
$\Z$-linear isomorphism $\eta^{(p)}: \Z^{2} \otimes_\Z \widehat{\Z}^{(p)} \iso \cT^{(p)}(A)$, where $\widehat{\Z}^{(p)}= \prod_{l\neq p} \Z_l$.\\
\\
The notation $\sim$ denotes up to a prime-to-$p$ isogeny.\\
\begin{thm}[Deligne--Rapoport]
 The functor $\cG^{(p)}$ is represented by a smooth pro-algebraic scheme $Sh^{(p)}_{/\Z_{(p)}}$. Moreover, 
there exists an isomorphism given by
$$ Sh^{(p)} \times \Q \iso Sh/G(\Z_p)_{/ \Q}.$$
(\cf\cite[\S 4.2.1]{Hi1}).\\
\end{thm}
\noindent \\
Let $g \in G(\A_{\Q,f}^{(p)})$ act on $Sh^{(p)}_{/\Z_{(p)}}$ via
\beq
x= (A,\eta^{(p)}) \mapsto gx= (A,\eta^{(p)}\circ g).
\eeq
\noindent\\
\\
For a sufficiently small $K$, we have the  corresponding level-$K^{(p)}$ Shimura variety $Sh^{(p)}_{K/\Z_{(p)}}$. 
For the moduli interpretation, we consider the functor $\cG_{K}^{(p)}$ essentially as above with  
the level $K^{(p)}$-structure instead of the hypothesis (PM2) (\cf\cite[\S4]{Hi1}).\\
\\
For geometrically irreducible schemes  
$Sh_{K/ \Z_{(p)}}^{(p)}(\mathfrak{c})$, we have 
\beq 
Sh_{K}^{(p)}= \bigsqcup_{ [\mathfrak{c}] \in \Cl_{\Q}^{+}(K)} Sh_{K}^{(p)} (\mathfrak{c}).
\eeq
Here $\Cl_{\Q}^{+}(K)$ is the narrow ray class group of $\Q$ with level $\det(K)$.\\
\\
\\
\subsubsection{Igusa tower}
\noindent In this subsection, we briefly recall the notion of 
$p$-ordinary Igusa tower over the $p$-integral modular Shimura variety. 
The Igusa tower underlies the geometric theory of $p$-adic modular forms. 
We refer to \cite[Ch. 8]{Hi1} for details.\\
\\
Let the notation and hypotheses be as in the introduction and \S2.1.2. 
In particular, $\cW$ is the strict Henselisation inside $\overline{\Q}$ of the local ring of $\Z_{(p)}$ corresponding to $\iota_p$, $W(\F)$ the Witt ring 
and $\F$ the residue field of $\cW$.\\
\\
Let $Sh^{(p)}_{/\cW} = Sh^{(p)}(G,X) \times_{\Z_{(p)}} \cW$ and $Sh^{(p)}_{/ \F}= Sh^{(p)}_{/\cW} \times_{\cW} \F$.\\
\\
From now, let $Sh$ (resp. $Sh_{K}$) denote $Sh^{(p)}_{/ W}$ or  $Sh^{(p)}_{/ \F}$ (resp. $Sh^{(p)}_{K/ W}$ or $Sh^{(p)}_{K/ \F}$). 
Unless otherwise stated, the base would be evident from the context. Let $\cA$ be the universal elliptic curve over $Sh$.\\
\\
Let $Sh^{ord}$ be the subscheme of $Sh$ on which the Hasse-invariant does not vanish. It is an open dense subscheme. 
Over $Sh^{ord}$, the connected part $\cA[p^m]^{\circ}$ of $\cA[p^m]$ is étale-locally isomorhpic to 
$\mu_{p^{m}}$.
We now define the Igusa tower. For $m \in \mathbb{N}$, the $m^{th}$-layer of the Igusa tower over $Sh^{ord}$ is defined by
\beq Ig_m=\Isom_{\Z_{p}}(\mu_{p^m}, \cA[p^m]^\circ).\eeq
Note that the projection $\pi_m: Ig_m \rightarrow Sh^{ord}$ is finite and étale. 
The full Igusa tower over $Sh^{ord}$ is defined by
\beq Ig=Ig_\infty=\varprojlim Ig_m = \Isom_{\Z_{p}}(\mu_{p^\infty}, 
\cA[p^\infty]^\circ).\eeq
($\acute{E}t$) Note that the projection $\pi: Ig \rightarrow Sh^{ord}$ is étale.\\
\\
Let V be an irreducible component of $Sh$ and $V^{ord}$ denote $V \cap Sh^{ord}$. 
Let $I$ be the inverse image of $V^{ord}$ under $\pi$. In \cite[Ch.8]{Hi1} and \cite{Hi2}, it has been shown that\\
\\
(Ir) $I$ is an irreducible component of $Ig$.\\
\\
For 
$\mathfrak{c}$
and $K$ as above, we can analogously define the $p$-ordinary Igusa tower $Ig_{K}(\mathfrak{c})$ over level-$K^{(p)}$ Shimura variety $Sh^{(p)}_{K}(\mathfrak{c})$.\\
\\
\\
\subsubsection{CM points}
In this subsection, we briefly recall some notation regarding CM points on the Igusa tower. 
The notion is originally due to Shimura and underlies the theory of Shimura varieties.
We refer to \cite{Sh3} for details.\\
\\
Let the notation and hypotheses be as in the introduction and \S2.1.3. 
In particular, $K/\Q$ is a $p$-ordinary CM quadratic extension. 
Let $\{1,\omega\}$ be a $\Z$-basis of $\cO_{K}$.
Let $\mathfrak{a} \subset \cO_{K}$ be a lattice. 
Let 
$$R(\mathfrak{a})=\bigg{\{}\alpha \in \cO_{K}\big{|} \alpha\mathfrak{a}\subset \mathfrak{a}\bigg{\}}$$ 
be the corresponding order of $K$ and $\mathfrak{f}(\mathfrak{a})\subset \Z$ the 
corresponding conductor ideal. Recall that $R(\mathfrak{a})=\Z+\mathfrak{f}(\mathfrak{a})\cO_{K}$.\\
\\
Let $\Sg$ be a $p$-ordinary CM type of K i.e. an embedding $\iota_{\infty}: K \hookrightarrow \C$. 
By CM theory of Shimura--Taniyama--Weil (\cf \cite{Sh3}), the complex torus $X(\mathfrak{a})(\C)=\C/\Sg(\mathfrak{a})$  
is algebraisable to a CM elliptic curve with CM type $(K,\Sg)$. 
Here $\Sg(\mathfrak{a})=\big{\{}\iota_{\infty}(a) |a \in \mathfrak{a}\big{\}}$. 
When the conductor $\mathfrak{f}(\mathfrak{a})$ is prime-to-$p$, the elliptic curve 
$X(\mathfrak{a})$ extends to an 
elliptic curve over $\cW$. We denote it by $X(\mathfrak{a})_{/\cW}$. 
In this case, a construction of a pair $(\eta^{(p)}(\mathfrak{a}),\eta_{p}^{ord}(\mathfrak{a}))_{/W}$ satisfying (PM2) and \S2.1.3 
is given in \cite[\S3.3]{Hs2}. 
By definition, this gives rise to a CM point
$$x(\mathfrak{a})=(X(\mathfrak{a}),\eta^{(p)}(\mathfrak{a}),\eta_{p}^{ord}(\mathfrak{a})) \in Ig(W).$$
\noindent\\
\\
\subsubsection{Tate objects} In this subsection, we briefly recall some notation
regarding Tate objects on the modular curves. 
They naturally arise from the construction of compactification of the modular curves. 
In turn, they give rise to the key notion of $q$-expansion of elliptic modular forms. 
We refer to \cite[\S 1.1]{Ka} and \cite[\S 4.1.5]{Hi1} for details.\\
\\ 
Let the notation and hypotheses be as in the introduction and \S2.1.3. 
Let $L$ be a lattice in $\Q$, $n$ a positive integer, 
$L_{n}=\big{\{}x\in L| x>-n\big{\}}$ and 
$A((L))=\varinjlim A\powerseries{L_{n}}$.
Pick two fractional ideals $\mathfrak{a},\mathfrak{b}$ of $\Q$ prime-to-$p$. 
To this pair, Tate associated a certain ellipitic curve  
$Tate_{\mathfrak{a},\mathfrak{b}}(q)_{/ \Z((\mathfrak{ab}))}$.  
Formally, $Tate_{\mathfrak{a},\mathfrak{b}}(q)=\mathfrak{a}^{*}\otimes_{\Z} \mathbb{G}_{m}/q^{\mathfrak{b}}$ for $\mathfrak{a}^{*}=\mathfrak{a}^{-1}$. 
It is endowed with a canonical prime-to-$p$ level structure $\eta_{can}^{(p)}$, 
$p^\infty$-level structure $\eta_{p,can}^{ord}$ and a generator $\omega_{can}$ of $\Omega_{Tate_{\mathfrak{a},\mathfrak{b}}(q)}$. 
We thus obtain a Tate object 
$$x(Tate_{\mathfrak{a,b}}(q))=(Tate_{\mathfrak{a},\mathfrak{b}}(q), \eta_{can}^{(p)}, \eta_{p,can}^{ord}).$$
\noindent\\
\\
\subsection{Modular forms}
\noindent In this section, we recall certain generalities about elliptic modular forms  
including the adelic notion.\\ 
\\
\subsubsection{Classical modular forms} 
In this subsection, we recall the notion of classical modular forms. 
We refer to \cite[\S4.2]{Hi1} and \cite[\S2.5]{Hs2} for details.\\
\\
\underline{Geometric definition}.
We first recall the geometric definition of classical modular forms.\\ 
\\
Let the notation and assumptions be as in \S2.1. 
Let $k$ be a positive integer and $R$ a $\Z_{(p)}$-algebra. Let $T$ be the torus $\mathbb{G}_{m,/\Z}$ and 
$\underline{k}$ the character of $T$ arising from 
$t \mapsto t^{k}$ for $t$ the co-ordinate.
A classical modular form of weight $k$ and level $K$ over $R$ is a 
function $f$ of isomorphism classes of $x=(\underline{\cA},\omega)$ where $\underline{\cA} \in Sh_{K}^{(p)}(S)$ and 
$\omega$ a differential form generating $H^{0}(A,\Omega_{A/S})$ over $S$ for an $R$-algebra $S$, locally under the Zariski topology such that the
following conditions are satisfied.\\
\\
{ (Gc1)} If $x \iso x'$, then $f(x)=f(x')\in S$.\\
{ (Gc2)} $f(x \otimes_{S} S') = \rho(f(x))$ for any $R$-algebra homomorphism $\rho: S \to S'$.\\
{ (Gc3)} $f(\underline{\cA},s\omega)=\underline{k}(s)^{-1}f(\underline{\cA},\omega)$, where $s \in T(S)$.\\
{ (Gc4)} $f|_{\mathfrak{a,b}}(q):=f(Tate_{\mathfrak{a},\mathfrak{b}}(q), \eta_{can}^{(p)},\omega_{can}) \in R\powerseries { q^{\mathfrak{ab}_{\geq0}}}$, 
where $\mathfrak{ab}_{\geq0}=(\mathfrak{ab}\cap \Q_{+}) \cup\{0\}$. We say that $f|_{\mathfrak{a,b}}(q)$ 
is the $q$-expansion of $f$ at the cusp $(\mathfrak{a,b})$.\\
\\
Let $M_{k}(K,R)$ be the space of the functions satisfying the above conditions (Gc1-4). 
For $ f \in M_{k}(K,R)$, we have the following fundamental $q$-expansion principle (\cf\cite[Thm. 4.21]{Hi1}).\\
\\
{ ($q$-exp)} The $q$-expansion map
$f \mapsto f|_{\mathfrak{a,b}}(q) \in R\powerseries { q^{\mathfrak{ab}_{\geq0}}}$ 
determines $f$ uniquely.\\
\\
\underline{Adelic definition}. 
We now briefly describe relation of the geometric definition of modular forms to the adelic counterpart.\\
\\
Let the notation and assumptions be as before. In particular, $k$ is a positive integer.
For $\tau \in X^{+}$ the upper half plane and $g=\MX{a}{b}{c}{d} \in G(\R)$, 
let 
$$
J(g, \tau)^{k} = (c\tau+d)^{k}.
$$
Let $f$ be a modular form of weight $k$ over $\C$. 
In view of the complex uniformisation of the modular curves in (2.1), 
we may regard $f$ as a function $f: X^{+} \times G(\A_{\Q,f}) \rightarrow \C$ (also see \cite[\S2.5.2]{Hs2}).
Let $\varphi: G(\A_{\Q,f}) \rightarrow \C$ be the function given by 
$$
f(\tau,g_{f}) = \varphi(g) \cdot J(g_{\infty}, i)^{k}(\det{g_{\infty}})^{-k}|\det{g}|_{\A}^{k/2}.
$$
Here 
$g_{\infty} \in G(\R)$ such that $g_{\infty}i=\tau$ and $\det{g_{\infty}}>0$. 
Moreover, $g=(g_{\infty},g_{f})$ and $|\cdot|_{\A}$ is the adelic norm. 
It turns out that $\varphi$ is an adelic modular form.\\
\\
Via the above identity, we may also begin with an adelic modular form and obtain a classical modular form.
Representation theoretic methods can often be used to construct adelic modular forms. For example, the construction of toric modular forms (\cf \cite{Hs3}) is based on adelic formalism (also see Appendix A).\\
\\
\\
\subsubsection{$p$-adic modular forms}
In this subsection, we recall the geometric definition of $p$-adic modular forms. 
The theory is due to Katz. We refer to \cite[\S8]{Hi1} for details.\\
\\
Let the notation and assumptions be as in \S2. 
Let $R$ be a $p$-adic algebra.  A $p$-adic modular form of level $K$ over $R$ is a function $f$ 
of isomorphism classes of $x=(\underline{\cA}, \eta_{p}^{ord}) \in Ig_{K}(S)$ defined over a  
$p$-adic $R$-algebra $S$ such that the following conditions are satisfied.\\
\\
{ (Gp1)} If $x \iso x'$, then $f(x)=f(x')\in S$.\\
{ (Gp2)} $f(x \otimes_{S} S') = \rho(f(x))$ for any $p$-adic $R$-algebra homomorphism $\rho: S \to S'$.\\
{ (Gp3)} $f|_{\mathfrak{a,b}}(q):=f(x(Tate_{\mathfrak{a,b}}(q))) \in R\powerseries { q^{\mathfrak{ab}_{\geq0}}}$. 
We say that $f|_{\mathfrak{a,b}}(q)$ 
is the $q$-expansion of $f$ at the cusp $(\mathfrak{a,b})$.\\
\\
 Let $V(K,R)$ be the space of the functions satisfying the above conditions (Gp1-3).  
We can regard a classical modular form as a $p$-adic modular form as follows. 
Let $f \in M_{k}(K,R)$
and $(\underline{\cA},\eta_{p}^{ord}) \in Ig_{K}(S)$. Note that $\eta_{p}^{ord}$ induces an isomorphism 
$\eta_{p}^{ord}: \Lie(\widehat{\mathbb{G}}_{m}) \to \Lie(A)$. 
It follows that $(\eta_{p}^{ord})^{*} (\frac{dt}{t})$ generates $H^{0}(A,\Omega_{A})$
as an $R$-module.
We regard 
$f$ as a $p$-adic modular form of level $K$ via 
$$f \mapsto \hat{f}(\underline{\cA},\eta_{p}^{ord}):= f(\underline{\cA},(\eta_{p}^{ord})^{*} (\frac{dt}{t})).$$
Moreover, $\hat{f}|_{\mathfrak{a,b}}(q)=f|_{\mathfrak{a,b}}(q)$ and we have an embedding 
$M_{k}(K,R)\hookrightarrow V(K,R)$.\\
\\
For $f \in V(K,R)$, we have the following fundamental $q$-expansion principle (\cf\cite[Thm. 4.21]{Hi1})\\
\\
{ ($q$-exp)'} The $q$-expansion map
$f \mapsto f|_{\mathfrak{a,b}}(q) \in R\powerseries { q^{\mathfrak{ab}_{\geq0}}}$ determines $f$ uniquely.\\
\\
\\
Another key fact is the completeness of the space $V(K,R)$ of $p$-adic modular forms under $q$-expansion topology.\\
\subsubsection{Nearly holomorphic modular forms} In this subsection, we recall generalities regarding nearly holomorphic forms. The theory is due to Shimura, Katz, Harris and Urban. We refer to \cite{U} and references therein for details.\\
\\
\underline{Geometric definition}.
Let the notation and assumptions be as in \S2.1. 
Let $k, r$ be non-negative integers and $R$ a $\Z_{(p)}$-algebra. 
Let $R[X]_{\leq r}$ be the submodule of polynomials over $R$ with degree at most $r$. 
Let $B$ be the Borel subgroup of $\SL_{2/\Z}$.  
Let $\underline{\rho^{r}_{k}}$ be the representation of $B(R)$ arising from the action on $R[X]_{\leq r}$ given by
$$
\MX{a}{b}{0}{a^{-1}} \cdot P(X) = a^{k} P(a^{-2}X + ba^{-1})
$$
for $P(X) \in R[X]_{\leq r}$.\\
\\
Let $S$ be an $R$-algebra. 
A nearly holomorphic modular form of weight $k$, order at most $r$ and level $K$ over $R$ is an $R[X]_{\leq r}$-valued  
function $f$ of isomorphism classes of $y=(\underline{\cA},\omega, \omega')$ where $\underline{\cA} \in Sh_{K}^{(p)}(S)$ and 
$(\omega,\omega')$ a basis of $H^{1}_{dR}(A_{/S})$. 
Here $\omega$ is a basis of $H^{0}(A,\Omega_{A/S})$ 
locally under the Zariski topology
with 
 $\langle \omega, \omega' \rangle_{dR} = 1$ 
 for the de Rham pairing $\langle \cdot, \cdot \rangle$
 such that the
following conditions are satisfied.\\
\\
{ (Gnh1)} If $y \iso y'$, then $f(y)=f(y')\in S$.\\
{ (Gnh2)} $f(y \otimes_{S} S') = \rho(f(y))$ for any $R$-algebra homomorphism $\rho: S \to S'$.\\
{ (Gnh3)} $f(\underline{\cA},s\cdot (\omega,\omega'))=\underline{\rho_{k}^{r}}(s)^{-1}\cdot f(\underline{\cA},\omega)$, where $s \in B(S)$.\\
{ (Gnh4)} $f|_{\mathfrak{a,b}}(q):=f(Tate_{\mathfrak{a},\mathfrak{b}}(q), \eta_{can}^{(p)},\omega_{can}, \omega'_{can}) \in R\powerseries { q^{\mathfrak{ab}_{\geq0}}}[X]_{\leq r}$, 
where $\omega'_{can}$ arises from $\omega_{\can}$ (\cf \S2.1.5) via the de Rham pairing as above and $\mathfrak{ab}_{\geq0}=(\mathfrak{ab}\cap \Q_{+}) \cup\{0\}$. We say that $f|_{\mathfrak{a,b}}(q)$ 
is the $q$-expansion of $f$ at the cusp $(\mathfrak{a,b})$.\\
\\
Let $N_{k}^{r}(K,R)$ be the space of the functions satisfying the above conditions (Gnh1-4). 
By definition, $N_{k}^{0}(K,R)$ is nothing but the space of classical modular forms $M_{k}(K,R)$.
For $ f \in N_{k}^{r}(K,R)$, we have the following fundamental $q$-expansion principle.\\
\\
{ ($q$-exp)} The $q$-expansion map
$f \mapsto f|_{\mathfrak{a,b}}(q) \in R\powerseries { q^{\mathfrak{ab}_{\geq0}}}[X]_{\leq r}$ 
determines $f$ uniquely.\\
\\
Let $$N_{k}(K,R)=\bigoplus_{r \geq 0} N_{k}^{r}(K,R)$$ be the space of nearly holomorphic modular forms with weight $k$, level $K$ over the base $R$.\\
\\
\underline{Differential operators} We briefly describe differential operators on the space of modular forms and its relation to nearly holomorphic modular forms.\\
\\
Let the notation and assumptions be as above. Let $\delta_{k}$ be the Maass--Shimura differential operator on the space of modular forms $N_{k}^{r}(K,\C)$ given by 
$$
\delta_{k}(f)=\frac{1}{2 \pi i} y^{-k} \frac{\partial}{\partial \tau} (y^{k}f)
$$
for $\tau=x + iy$ the complex variable. Note that $\delta_{k}(f) \in N_{k+2}^{r+1}(K,\C)$ with $X$ being $\frac{-1}{4 \pi y}$. For a positive integer $s$, let 
$$
\delta_{k}^{s}= \delta_{k+2s-2} \circ \cdot \cdot \cdot \circ \delta_{k}.
$$
We accordingly have $\delta_{k}^{s}(f) \in N_{k+2s}^{r+s} (K,\C)$.
The differential operator $\delta_{k}$ has a geometric interpretation and can in fact be extended $p$-integrally to $N_{k}^{r}(K,R)$ (\cf \cite [\S2.4]{U}). In particular, it can be defined on $f \in M_{k}(K,R)$ and we have  $\delta_{k}^{s}(f) \in N_{k+2s}^{r+s} (K,R)$.\\
\\
When $f$ is a newform and $\pi$ is the cuspidal automorphic representation of $\GL_2(\A)$ generated by $f$ (\S2.2.1), we indeed have $\delta_{k}^{s}(f) \in \pi$. \\
\\
On the other hand, we have Katz $p$-adic differential operator $d$ on the space of $p$-adic modular forms $V(K,R)$ with action on the $q$-expansion given by 
$$
d(g)|_{\mathfrak{a,b}}(q)=q\frac{d}{dq} (g|_{\mathfrak{a,b}}(q)).
$$
for $g \in V(K,R)$ (\cf \cite[1.3.6]{Hi5}).\\
\\
In view of the $q$-expansion principle, we have 
$$
d(f) = \delta_{k}(f)
$$
for $f \in M_{k}(K,R)$. We may thus regard the nearly holomorphic form $\delta_{k}^{s}(f)$ as an element in the space of $p$-adic modular forms $V(K,R)$. In this sense, we may regard the nearly holomorphic form $\delta_{k}^{s}(f)$ as being scalar valued.\\
\\
\subsubsection{Hecke operators}
In this subsection, we recall the definition of certain Hecke operators on the space of $p$-adic modular forms.\\
\\
Let the notation and assumptions be as in \S2. 
For $K$ as in \S2.2 and $g \in G(\A_{\Q,f}^{(p)})$, let $^{g}K=gKg^{-1}$. Note that the action $x \mapsto gx$ (\cf(2.3)) gives rise to an isomorphism $Sh_{K} \iso Sh_{^{g}K}$. 
For $f \in V(K,R)$, let $f|g$ be given by
\beq
(f|g)(x)=f(gx). 
\eeq
It follows that $f|g \in V(^{g}K,R)$.\\
\\
Recall $x=(\underline{\cA}, \eta_{p}^{ord}) \in Ig_{K}(S)$ (\S2.2.2). Let $(x_{i})_{0 \leq i \leq p}$ denote the points on the Igusa tower arising from distinct $p$-isogenies on the underlying elliptic curve $E$ with $x_{0}$ arising from the isogeny with kernel being the canonical subgroup of $E_{/S}$. Let 
$$
Ux= \frac{1}{p} \cdot \sum_{i=1}^{p} x_{i}, Vx= \big{(}E_{0},\frac{1}{p}\cdot\eta^{(p)},p\cdot\eta_{p}^{ord}\big{)}.
$$
\\
Let $f|U$ and $f|V$ be given by
$$
(f|U)(x)=f(Ux), (f|V)(x) = f(Vx).
$$
We have $f|U, f|V \in V(K,R)$ (\cite[\S3.8]{BDP1}).
\\
\section{Anticyclotomic toric periods} 
\noindent In this section, we consider non-triviality of anticyclotomic toric periods.
In \S3.1, we firstly describe the notion of anticyclotomic toric periods associated to a class of modular forms.  
We then discuss their non-triviality due to Hida. 
We provide a slightly different formulation and mild strengthening.
In \S3.2, we describe a linear independence of mod $p$ modular forms arising from Chai's theory of Hecke-stable subvarieties of a mod $p$ Shimura variety. The non-triviality is fundamentally based on the independence.\\
\\
\subsection{Non-triviality}
\noindent In this subsection, we consider non-triviality of anticyclotomic toric periods. For a related consideration, we also refer to \cite{Bu2}.\\
\\
Let the notation and hypotheses be as in \S1-2. 
In particular, $p$ is an odd prime, $N$ a positive integer with $p \nmid N$ and $c$ a positive integer with $(c, pN)=1$. Moreover, $K$ is an imaginary quadratic field with $p$ split, the embeddings $\iota_{\infty}\colon \Qbar \to \C$ and $\iota_{p}\colon \Qbar \to \C_p$. 
Further, $K_{\ell^n}^{-}$ is the anticyclotomic extension of $K$ with conductor $cN \cdot \ell^n$ for a prime $\ell \nmid 6pNc$, 
$\Gamma_{n}^{-}=\Gal(K_{\ell^n}^{-}/K)$ the anticyclotomic Galois group 
and $\Gamma_{\ell}^{-}=\varprojlim_{n} \Gamma_{n}^-$.\\
\\
Let $R_{n}=\Z+ cN\cdot\ell^{n}\cO_{K}$. Note that $\Pic(R_n)=\Gamma_{n}^-$. 
Let $$U_{n}= \C_{1} \times (\widehat{R_{n}}) \subset \C^{\times} \times {\bf{A}}_{K,f}^{\times}$$ be a compact 
subgroup for $\C_1$ the complex unit circle. 
Via the reciprocity map (\cf Notation), we identify $$\Gamma_{n}^{-} = K^{\times}{\bf{A}}^{\times}\backslash {\bf{A}}_{K}^{\times} / U_{n}.$$
Let $[\cdot]_{n}: {\bf{A}}_{K}^{\times} \rightarrow \Gamma_{n}^-$ be the quotient map. 
For $a \in {\bf{A}}_{K}^{\times}$, let $x_{n}(a)$ be the CM point $x([a]_n )$ associated to the ideal class $[a]_n$ as defined in \S2.1.4. 
Let ${\bf{c}}(a)$ be polarisation ideal of the CM point $x_{0}(a)$ defined in \cite[\S3.4]{Hs2}.
Also, let $[a]=\varprojlim_{n} [a]_{n} \in \Gamma_{\ell}^{-}$.\\
\\ 
Let $\lam$ be an arithmetic Hecke character over $K$ of infinity type $(j+\kappa_{c},-\kappa_{c})$ with respect to the embedding $\iota_\infty$ for $j, \kappa_{c} \geq 0$.
Let $\widehat{\lam}: {\bf{A}}_{K,f}^{\times}/K^{\times} \rightarrow \C_{p}^\times$ 
be its $p$-adic avatar given by 
$$  \widehat{\lam}(x)=\lam(x) x_{p}^{(j+\kappa_{c}(1-c))}                  .$$
Here $x \in {\bf{A}}_{K,f}^{\times}$, $x_p$ its $p$-component and $c \in \Gal(K/\Q)$ the non-trivial element.\\
\\
Let $f \in V(K_{0}(\ell), \cO)$ be a $p$-adic modular form for a finite flat extension $\cO$ over $\Z_p$ 
(\cf \S2.2.2) 
such that the following condition is satisfied:\\
\\
(MC) $f(x_{n}(ta))= \widehat{\lam}(a)^{-1}f(x_n(t))$ for $a \in U_{n}{\bf{A}}_{\Q,f}^\times$.\\
\\
For $\phi:\Gamma_{n}^{-} \rightarrow \overline{\Z}_p$, let $P_{f,\lambda}(\phi,n)$ be the toric period given by 
\beq
P_{f,\lambda}(\phi,n) = \sum_{[t]_{n} \in \Gamma_{n}^-} f(x_{n}(t))\widehat{\lam}(t)\phi([t]_{n}) .
\eeq
In view of the condition (MC), the expression on the right hand side is well defined if we change $t$ to $ta$, for $a \in U_{n} \Q^\times$. 
The toric period typically depends on $n$ and in general, the above definition may not give rise to a measure on $\Gamma_{\ell}^-$. 
Under additional hypothesis on $f$, we can normalise the toric periods and thereby obtain a measure on $\Gamma_{\ell}^-$ (\cf\cite[\S3.1 and \S3.4]{Hip}).\\
\\
Let $\mathfrak{X}_{\ell}$ denote the set of $\ell$-power order anticyclotomic characters of $\Gamma_\ell^-$. 
As $\nu \in \mathfrak{X}_{\ell}$ varies, we consider the non-triviality of toric periods $P_{f,\lam}(\nu,n)$ modulo $p$.  
For the rest of the section, we suppose that $f$ is a $p$-integral nearly holomorphic modular form defined over a number field (\cf \S2.2.3).\\
\\
To discuss non-triviality of the toric periods, we introduce further notation. 
Let $\Delta_{\ell}$ be the torsion subgroup of $\Gamma_{\ell}^-$, 
$\Gamma_{\ell}^{alg}$ the subgroup of $\Gamma_{\ell}^-$ generated by $[a]$ for $a \in ({\bf{A}}_{K}^{(\ell p)})^\times$ and 
$\Delta_{\ell}^{alg}= \Gamma_{\ell}^{alg} \cap \Delta_{\ell}$. 
Let $\mathcal{B}$ be a set of representatives of $\Delta_{\ell}/\Delta_{\ell}^{alg}$ 
and $\mathcal{R}$ a set of representatives of $\Delta_{\ell}^{alg}$ in $({\bf{A}}_{K}^{(\ell p)})^\times$. 
Let $\rho: {\bf{A}}_{K}^\times \hookrightarrow \GL_{2}({\bf{A}})$ be a torus embedding and 
$\varsigma \in \GL_{2}({\bf{A}}_{K}^{(\ell p)})$ as defined in \cite[\S3.1 and \S3.2]{Hs2}, respectively. 
For $a \in ({\bf{A}}_{K}^{(\ell p)})^\times$, let
\beq
f|[a]=f|\rho_{\varsigma}(a), \rho_{\varsigma}(a)= \varsigma^{-1}\rho(a) \varsigma 
\eeq
(\cf 2.2.4). 
Finally, let
\beq
f^{\mathcal{R}}= \sum_{r \in \mathcal{R}} \widehat{\lam}(r)f|[r]. 
\eeq
We consider the following hypothesis.\\
\\
(H) Let $\ell^{r}$ be the order of the $\ell$-Sylow subgroup of $\F_{p}[g,\lambda]^\times$ with $\F_{p}[g,\lambda]$ the mod $p$ Hecke field corresponding to the pair $(g,\lam)$.
For every $u \in \Z$ prime to $\ell$, there exists a positive integer $\beta \equiv u \mod{\ell^{r}}$ 
and $a \in  {\bf{A}}_{K,f}^{\times}$
such that 
$$\bfa_{\beta}(f^{\mathcal{R}},{\bf{c}}(a)) \neq 0 \mod{\mathfrak{m}_{p}}.$$ 
Here the Fourier expansion is at the cusp $(\Z,{\bf{c}}(a))$ (\S 2.2.3).\\
\\
We have the following fundamental result regarding non-triviality of the toric periods modulo $p$.\\
\begin{thm}(Hida) 
Let $p$ be an odd prime split in an imaginary quadratic field $K$. 
Let $\ell$ be another prime and 
$\mathfrak{X}_{\ell}$ the set of $\ell$-power order anticyclotomic character of the Galois group 
$\Gamma_{\ell}^-$ as above. 
Let $f$ be an algebraic $p$-integral nearly holomorphic modular form and $\lam$ an arithmetic Hecke character over $K$ as above such that (MC) holds. 
For $\nu \in \mathfrak{X}_\ell$ factoring though $\Gamma_n^-$, let $P_{f,\lam}(\nu,n)$ be the anticyclotomic toric period as above.

Suppose that the hypothesis (H) holds. Then, we have\\
$$P_{f,\lam}(\nu,n)  \neq 0 \mod{\mathfrak{m}_{p}},$$ for 
all but finitely many $\nu \in \mathfrak{X}_{\ell}$ factoring through $\Gamma_n^-$ as $n \rightarrow \infty$ (\cf\cite[Thm. 3.2]{Hip}).\\
\end{thm}
\begin{proof} 
This is essentially proven in \cite[Thm. 3.2]{Hip}. We outline a sketch.\\
\\
Let $\nu \in \mathfrak{X}_{\ell}$ minimally factor through $\Gamma_{n(\nu)}^{-}$.\\
\\
We fix a decomposition
$$
\Gamma_{\ell}^{-}=\Delta_{\ell} \times \Gamma
$$
with $\Gamma\simeq \Z_{\ell}$. Let $\Gamma_{n}$ be the image of $\Gamma$ in $\Gamma_{n}^{-}$ via 
$\Gamma \twoheadrightarrow \Gamma_{n}^{-}$. For a sufficiently large $n$, $\Delta_{\ell} \hookrightarrow \Gamma_{n}^{-}$ and 
\begin{equation}
\Gamma_{n}^{-}=\Delta_{\ell} \times \Gamma_{n}.
\end{equation}
We choose representatives $\mathfrak{b} \in \mathcal{B}$
\footnote{In the proof, we consider ideal-theoretic formulation in contrast with the earlier adelic formulation and accordingly regard Hecke characters as ideal class characters. 
}
prime to $p\ell N$ such that the norms $N(\mathfrak{b})$ are all different.\\
\\
Let $\F_{p^{r}}$ be the residue field of $E_{f}H$ at the prime above $p$ determined via the embedding 
$\iota_{p}$. From Shimura's reciprocity law, 
\beq
\bar{P}_{f,\lam}(\nu,n)  \neq 0 \iff \bar{P}_{f,\lam}(\nu^{\sigma},n) \neq 0
\eeq
for $\sigma \in \Gal(\F/\F_{p^{r}})$ (\cite[(3.10)]{Hip}). 
Here and in the remaining proof, we let $\bar{P}_{f,\lam}(\nu,n)$ denote 
$P_{f,\lam}(\nu,n)  \mod{\mathfrak{m}_{p}}$.
As for the reciprocity law, we rely on the hypothesis that 
the nearly holomorphic modular form $f$ is defined over a number field. \\
\\
Suppose that the non-vanishing does not hold for an infinite subset $\mathfrak{X}_{\ell}' \subset \mathfrak{X}_{\ell}$. We may suppose that for each $\nu \in \mathfrak{X}_{\ell}'$ the decomposition (3.4) holds with $n=n(\nu)$. 
In view of (3.4) and (3.5), we note
\begin{equation}
\sum_{\mathfrak{b} \in \mathcal{B}} \widehat{\lambda}(\mathfrak{b}) \sum_{\mathfrak{g}_{n} \in \Gamma_{n},\nu(\mathfrak{g}_{n}) \in \mu_{\ell^{r}}} \nu(\mathfrak{g}_{n}^{-1}) \cdot (\bar{f}^{\mathcal{R}}|[N(\mathfrak{b})])(x(\mathfrak{b}_{n}\mathfrak{g}_{n}\mathfrak{a})_{/\F}) = 0
\end{equation}
for $\mathfrak{b}_{n}$ the projection of $\mathfrak{b}$ to $\Gamma_{n}$ and any fixed $\mathfrak{a} \in \Gamma_{n}$(\cite[(3.15)]{Hip}). Moreover, $\bar{f}^{\mathcal{R}}$ denotes the mod $p$ reduction of $f^{\mathcal{R}}$ and 
$x(\cdot)_{/\F}$ the mod $p$ reduction of CM point $x(\cdot)$.\\
\\
For $n=n(\nu)$ sufficiently large, 
$$
\nu(\mathfrak{g}_{n}) \in \mu_{\ell^{r}} \iff \mathfrak{g}_{n} \in \ker (\Gamma_{n} \rightarrow \Gamma_{n-r}) =: G_{n,r}.
$$
Replacing $\mathfrak{X}_{\ell}'$ by possibly complement of a finite subset, we suppose that the above holds for each $\nu \in \mathfrak{X}_{\ell}'$ and also that 
$$
n \geq 2r.
$$
From the last inequality, we have an isomorphism $\Z/\ell^{n} \simeq G_{n,r}$ arising from 
$$
i \mapsto e_{n-r}\cdot R_{n}, e_{n,i}:=1+i \cdot cN\ell^{n} \cdot \omega.
$$
Replacing $\mathfrak{X}_{\ell}'$ by possibly an infinite subset, we suppose that the primitive $\ell^{r}$-th root of unity $\nu(e_{n-r}\cdot R_{n})$ is independent of $\nu$. Let $\zeta$ denote this primitive $\ell^{r}$-th root of unity and 
$$
n_{1}=\min_{\nu \in \mathfrak{X}_{\ell}'} (n(\nu)-r).
$$
In view of the preceding paragraph and (3.6), we note 
\begin{equation}
\sum_{\mathfrak{b} \in \mathcal{B}} \widehat{\lambda}(\mathfrak{b}) \cdot 
(\bar{f}^{\mathcal{R}}|N(\mathfrak{b}))(x(\mathfrak{b}_{n} \cdot \mathfrak{a})_{/\F}) = 0.
\end{equation}
 Here 
$$
\bar{f}^{\mathcal{R}} = \sum_{i=0}^{\ell^{r}-1} \zeta^{-i} \bar{f}^{\mathcal{R}}|B_{r,i}
\text{,  } 
B_{r,i}=\MX{1}{-i/\ell^{r}}{0}{1}
$$
and $\mathfrak{a} \in \ker(\Gamma_{n} \rightarrow \Gamma_{n_{1}})=\ker(\Gamma_{n}^{-} \rightarrow \Gamma_{n_{1}}^{-})$ a fixed element
(\cite[Proof of Thm. 3.2]{Hip}). 

From the hypothesis (H), the mod $p$ modular form $\bar{f}^{\mathcal{R}}$ is non-constant
(\cite[Proof of Thm. 3.2]{Hip}). 
Thus, the linear dependence in (3.7) contradicts Theorem 3.2.\\
\end{proof}
\noindent\\
\begin{remark} 
(1). In \cite{Hip}, Hida constructs a modular measure interpolating the normalised toric periods 
and discusses the non-triviality of the measure modulo $p$. 
The hypothesis that $g$ is an $U_{\ell}$-eigenform with 
the corresponding eigenvalue being a $p$-unit is only needed in the construction of the modular measure. 
For the sake of non-triviality of the toric periods alone, the hypothesis is not necessary.  
The author is grateful to Hida for this key remark.\\
\\
(2). As the theorem does not require the underlying nearly holomorphic form to be an $U_{\ell}$-eigenform, it can be used to slightly simplify the proof of $(\ell,p)$-non-triviality of 
central L-values in \cite{Hs2} and \cite{Hs3}.\\
\end{remark}
\noindent\\
\subsection{Linear independence}
\noindent In this subsection, we describe a linear independence of mod $p$ modular forms arising from Chai's theory of Hecke-stable subvarieties of a mod $p$ Shimura variety. The non-triviality of anticyclotomic toric periods is fundamentally based on the independence.\\
\\
Let the notation and hypotheses be as in \S3.1 \footnote{prior to the proof of Theorem 3.1}. Let $V=V_{/\F}$ be an irreducible component of the Igusa tower $Ig=Ig_{/\F}$ (\S2.1.3). Let $\underline{n}=0<n_{1}<n_{2}<...$ be an infinite sequence of non-negative integers and 
\begin{equation}
\Xi=\bigg{\{}
x(\sigma)_{/\F} \big{|} \sigma \in \Gamma_{n_{j}}^{-} \text{ with $\sigma \in \ker(\Gamma_{n_{j}}^{-} \rightarrow \Gamma_{r}^{-}), j \in \Z_{\geq 1}$} 
\bigg{\}}
\end{equation}
an infinite set of mod $p$ CM points for an integer $r$ with $0 \leq r \leq n_{1}$.\\
\\
Let $\mathcal{C}_{\Xi}$ be the space of functions over $\Xi$ with values in $\mathbb{P}^{1}(\F)=\F \cup \{\infty\}$.
In view of the Zariski density of $\Xi$ in $V$, we have an embedding $$\F(V) \hookrightarrow\mathcal{C}_{\Xi}$$
for $\F(V)$ the function field.
We let the class group $\Gamma_{\ell}^{-}$ act on $\mathcal{C}_{\Xi}$ via 
$$
\sigma' \cdot g(x(\sigma)_{/\F}) = g(x(\sigma \cdot \sigma')_{/\F}).
$$
We have the following linear independence.\\
\begin{thm}(Hida)
Let the notation be as above. 
Let $\mathcal{L}$ be a line bundle over $Ig$. 
Let $\Delta \subset \Gamma_{\ell}^{-}$ be a finite subset independent modulo 
$\Gamma_{\ell}^{alg} \subset \Gamma_{\ell}^{-}$
and $\{ s_{\delta} \in \mathcal{L} \}_{\delta \in \Delta}$ a set of non-constant global sections $s_{\delta}$ finite at 
$\Xi$. Then,
$$
\{\delta \cdot s_{\delta}\}_{\delta \in \Delta} \subset \mathcal{C}_{\Xi}
$$
is linearly independent (\cite[Cor. 2.9]{Hip}).\\
\end{thm}
\begin{remark}
The Chai--Oort rigidity principle predicts that a Hecke stable subvariety of a mod $p$ Shimura variety is a Shimura subvariety. 
The principle has been studied by Chai in a series of articles (\cf \cite{Ch1}, \cite{Ch2} and \cite{Ch3}). The theorem is a consequence of a variant of the principle for self-products of mod $p$ modular curves where Hecke action is replaced with its local analogue (\cite[Prop. 2.8]{Hip}).\\
\end{remark}
\noindent\\
\section{Non-triviality of $p$-adic Abel--Jacobi image modulo $p$}
\noindent In this section, we prove the non-triviality of generalised Heegner cycles modulo $p$. 
For simplicity, we first consider the case of Heegner points. 
In \S4.1, we first state the $p$-adic Waldspurger formula  
and then prove the non-triviality of the $p$-adic formal group logarithm of Heegner points modulo $p$ (\cf Theorem A). 
In \S4.2, we describe generalities regarding generalised Heegner cycles and state conjectures regarding their non-triviality.
In \S4.3, we prove non-triviality of the $p$-adic Abel--Jacobi image of generalised Heegner cycles modulo $p$.\\
\\
\subsection{Non-triviality, I} 
In this subsection,  we first state the $p$-adic Waldspurger formula due to Bertolini--Darmon--Prasanna 
and then prove the non-triviality of the $p$-adic formal group logarithm of Heegner points modulo $p$ (\cf Theorem A).\\
\\
Let the notation and hypotheses be as in the introduction. 
In particular, $p$ is an odd prime, $N$ a positive integer with $p \nmid N$, $b$ a positive integer with $(b, pN)=1$ and $\ell$ a prime such that $\ell \nmid 2pNb$. Moreover, $K$ is an imaginary quadratic field with $p$ split satisfying the Heegner hypothesis (Hg), 
$K_{\ell^n}^{-}$ is the anticyclotomic extension of $K$ with conductor $bN \cdot \ell^n$, 
$\Gamma_{n}^{-}=\Gal(K_{\ell^n}^{-}/K)$ the anticyclotomic Galois group and $\Gamma_{\ell}^{-}=\varprojlim_{n} \Gamma_{n}^-$.\\
\\
Further, $f \in S_{2}(\Gamma_{0}(N),\epsilon)$ is a weight two Hecke-eigen cuspform with $T_{p}$-eigenvalue $\bfa_{p}(f)$ and $B_{f}$ an abelian variety associated to $f$. 
Also, $\chi$ is a finite order Hecke character over $K$ such that $\chi {\bf{N}}_{K} \in \Sg_{cc}(b,\mathfrak{N},\epsilon)$ for ${\bf{N}}_K$ the norm Hecke character over $K$ and $H_{\chi}$ the abelian extension of $K$ cut out by $\chi$. 
Let $E_{f,\chi}$ be the Hecke field corresponding to the pair $(f,\chi)$. 
The restriction of the $p$-adic logarithm gives a homomorphism 
$\log_{\omega_{B_{f}}}:  B_{f}(H_{\chi}) \rightarrow \C_{p}$ for the invariant differential $\omega_{B_{f}}$ arising from the newform $f$ and 
we extend it to 
$$
\log_{\omega_{B_{f}}}: B_{f}(H_{\chi}) \otimes_{\Z_{f}} E_{f,\chi} \rightarrow \C_{p}
$$ 
via $E_{f,\chi}$-linearity.\\
\\
For $\sigma \in \Gal(H_{bN}/K)$, let $x(\sigma) \in Ig_{\Gamma_{0}(N)}(W)$ be the CM point as in \S2.1.4.\\
\\
\\
We recall $\epsilon(f, \chi^{-1})=-1$ for the global root number $\epsilon(f, \chi^{-1})$ corresponding to Rankin--Selberg convolution $f \times \chi^{-1}{\bf{N}}_{K}^{-1}$ of the pair $(f,\chi^{-1})$. 
The generalised Heegner hypothesis (Hg) holds. The image of CM points under modular parametrisation of $B_{f}$ gives rise to a Heegner point 
$$
P_{f}(\chi)\in B_{f}(H_{\chi}) \otimes_{\Z_{f}} E_{f,\chi}.
$$
For the Hecke operators $V$ and $U$ on the space $V(\Gamma_{0}(N),R)$ as in \S 2.2.4  
and $g\in V(\Gamma_{0}(N),R) $, let $g^{(p)}$ be the $p$-depletion given by 
\beq
g^{(p)}=g|(VU-UV).
\eeq
The Fourier expansion of $g^{(p)}$ at a cusp ${\bf{c}}$ is given by
\beq
g^{(p)}(q)=\sum_{p\ndivide n} \bfa_{n}(g, {\bf{c}})q^{n} 
\eeq
(\cf \cite[(3.8.4)]{BDP1}).\\
\\
Let $d$ be the Katz $p$-adic differential operator as in \S2.2.3. 
In view of (4.2) and completeness of the space $V(\Gamma_{0}(N),R)$ (\cf \S2.2.2), it follows that $$d^{-1}(g^{(p)})=\lim_{j \rightarrow 0} d^{-1+j}(g^{(p)})$$ is a weight zero $p$-adic modular form. Here the limit is taken with respect to the $p$-adic topology. Indeed, limit of $q$-expansions of 
the nearly holomorphic forms $d^{-1+j}(g^{(p)})$ at the cusp exists with respect to the $p$-adic topology.  
The existence can be readily seen from action of the $p$-adic differential operator on the $q$-expansions (\S2.2.3). 
In view of the completeness, it then follows that the limit in fact arises from $q$-expansion of a $p$-adic modular form.\\
\\
We have the following $p$-adic Waldspurger formula for anticyclotomic toric periods of the weight zero $p$-adic modular form $d^{-1}(f^{(p)})$.\\
\begin{thm} (Bertolini--Darmon--Prasanna) 
Let $p$ be an odd prime and $N$ a positive integer such that $p \ndivide N$. 
Let $K$ be an imaginary quadratic field satisfying (ord) and (Hg) with $\mathfrak{p}$ a prime above $p$ determined via an initial embedding $\iota_p: \overline{\Q} \hookrightarrow \C_p$.
Let $f\in S_{2}(\Gamma_{0}(N),\epsilon)$ be a normalised newform with $T_{p}$-eigenvalue $\bfa_{p}(f)$ and $\eta$ a finite order Hecke character such that $\eta {\bf{N}}_{K} \in \Sg_{cc}(b,\mathfrak{N},\epsilon)$ 
for some $b$ prime to $pN$. We then have
$$ \sum_{\sg \in \Gal(H_{bN}/K)} \chi(\sg) d^{-1}(f^{(p)})( x(\sg)) = (1-\chi^{-1}(\overline{\mathfrak{p}})p^{-1}\bfa_{p}(f)+\chi^{-2}(\overline{\mathfrak{p}})\epsilon(p)p^{-1}) \cdot \log_{\omega_{B_{f}}}(P_{f}(\chi)).$$
\noindent\\
\end{thm}
\begin{proof} This is essentially proven in \cite{BDP1} and \cite{BDP2}. 
We give some details as the result in \cite{BDP1} and \cite{BDP2} explicitly involves a Rankin--Selberg $p$-adic L-function.\\
\\
We first recall that 
\beq
T_{p}(f)=\bfa_{p}(f)f, \langle p \rangle (f)=\epsilon(p)f.
\eeq
\noindent\\
Here $\langle p \rangle$ denotes the diamond operator associated with $p$.\\
\\
Let $F_{f}$ be the Coleman primitive of the differential $\omega_{f}$ associated with $f$ (\cf \cite[\S3.7]{BDP1}). 
From Coleman integration (\cf \cite[Lem. 3.23]{BDP1} and \cite[Thm. 3.12]{BDP2}), we have
\beq
F_{f}(x(\sg))=\log_{\omega_{B_{f}}}(\Phi_{f}(\Delta_{\sg})).
\eeq 
Here $\Delta_{\sigma}=x(\sg) - (\infty) \in J_{1}(N)(H_{bN})$ under the embedding $X_{1}(N)\hookrightarrow J_{1}(N)$.\\
\\
On the ordinary locus, the Coleman primitive $F_{f}$ satisfies
$$
dF_{f}=f.
$$
In view of (4.3) and the definition of $p$-depletion (\cf \cite[Lem. 3.23 and Prop. 3.24]{BDP1}), we have 
\beq
F_{f}^{(p)}(x(\sg))= F_{f}( x(\sg)) - \epsilon(p)\bfa_{p}(f) F_{f}(\mathfrak{p}*x(\sg))+\frac{\epsilon(p)}{p}F_{f}(\mathfrak{p}^{2}*x(\sg)).
\eeq 
Here the $``*"$-action of $\mathfrak{a}$ on the CM point $x(\sigma)$ is given by the triple arising from CM elliptic curve
$\frac{X(\sigma)} {X(\sigma)[\mathfrak{a}]}$ for $\mathfrak{a=p,p}^{2}$ (\cf \cite[(1.4.8)]{BDP1}). As $\mathfrak{a}$ is coprime to $N$, the level $\Gamma_{0}(N)$-structure on $X(\sigma)$ gives rise to level $\Gamma_{0}(N)$-structure on $\frac{X(\sigma)}{X(\sigma)[\mathfrak{a}]}$.\\
\\
Along with (4.4), the summation of the above identity over $\sg \in \Gal(H_{bN}/K)$ finishes the proof 
(compare \cite[Proof of Thm. 5.13]{BDP1}
\footnote{As opposed to \cite[Thm. 5.13]{BDP1}, $\chi(\sigma)$ appears in the left hand side of the formula in Theorem 4.1 due to geometrically normalised reciprocity law (\cf Notation).}).\\
\end{proof}
\noindent\\
We now consider the variation of the above identity over the $\Z_\ell$-anticyclotomic extension of $K$ for a prime 
$\ell \nmid 2pN$ as above. 
We fix a finite order Hecke character $\eta$ such that $\eta {\bf{N}}_{K} \in \Sg_{cc}(c,\mathfrak{N},\epsilon)$ for some $c$ satisfying  $(c,pN)=1$.
In the remaining subsection, we apply the formalism of anticyclotomic toric periods (\cf \S3) to the pair $(d^{-1}(f^{(p)}),\eta)$.\\
\\
Recall that $\mathfrak{X}_{\ell}$ denotes the set of $\ell$-power order anticyclotomic characters of $\Gamma_{\ell}^-$. 
As $\nu \in \mathfrak{X}_{\ell}$ varies, the global root number $\epsilon(f,\eta^{-1}\nu^{-1})$ still equals $-1$ (Lemma 4.4). \\
\\
The $p$-adic modular form $d^{-1}(f^{(p)})$ being of weight zero, 
the hypothesis (MC) holds and the toric periods are thus well defined.
The Hecke character $\eta$ being of finite order with conductor $c\mathfrak{N}_\epsilon$, for $\phi:\Gamma_{n}^{-} \rightarrow \overline{\Z}_p$ we have
\beq
P_{d^{-1}(f^{(p)}),\eta}(\phi,n)= \sum_{[t]_{n} \in \Gamma_{n}^{-}} (\phi\eta)([t]_{n}) d^{-1}(f^{(p)})(x_{n}(t)).
\eeq
The following is an $(\ell,p)$-family version of the $p$-adic Waldspurger formula.\\
\begin{cor} Let the notation and hypotheses be as above. For $\nu \in \mathfrak{X}_{\ell}$, we have 
\beq
P_{d^{-1}(f^{(p)}),\eta}(\nu,n) = (1-(\eta\nu)^{-1}(\overline{\mathfrak{p}})p^{-1}\bfa_{p}(f)+(\eta\nu)^{-2}(\overline{\mathfrak{p}})\epsilon(p)p^{-1})\log_{\omega_{B_{f}}}(P_{f}(\eta\nu)).
\eeq
\end{cor}
\begin{proof} 
We first note that as $\eta {\bf{N}}_{K} \in \Sg_{cc}(c,\mathfrak{N},\epsilon)$, it follows that $\eta\nu {\bf{N}}_{K} \in \Sg_{cc}(c\ell^{k},\mathfrak{N},\epsilon)$ 
for a non-negative integer $k$. 
The lemma thus immediately follows from the previous theorem.\\
\end{proof}
\noindent\\
We now consider non-triviality of $P_{d^{-1}(f^{(p)}),\eta}(\nu,n)$ modulo $p$. 
We are now ready to prove Theorem A regarding $p$-indivisibility of Heegner points (\cf \S1).\\
\begin{thm}
Let $p$ be an odd prime and $N$ a positive integer such that $p \ndivide N$. 
Let $K$ be an imaginary quadratic field satisfying (ord) and (Hg) with $\mathfrak{p}$ a prime above $p$ determined via an initial embedding $\iota_p: \overline{\Q} \hookrightarrow \C_p$. 
Let $f\in S_{2}(\Gamma_{0}(N),\epsilon)$ be a normalised newform with $T_{p}$-eigenvalue $\bfa_{p}(f)$ and $\eta$ a finite order Hecke character such that $\eta {\bf{N}}_{K} \in \Sg_{cc}(c,\mathfrak{N},\epsilon)$ 
for some $c$ prime to $pN$.
Suppose that\\
\\
(1). The residual representation $\rho_{f}|_{G_{K}} \mod{\mathfrak{m}_{p}}$ is absolutely irreducible,\\
(2). $N^-$ is square-free and\\
(3). $p \ndivide \prod_{v|c^{-}}|\Delta_{\eta,v}|$ for $\Delta_{\eta,v}$ the finite group $\eta(\cO_{K_{v}}^{\times})$.\\
\\
Let $\ell\neq p$ be an odd prime unramified in $K$ and prime to $cN$. 
Let $\mathfrak{X}_{\ell}$ be the set of $\ell$-power order anticyclotomic Hecke characters over $K$ as above. 
Then, for all but finitely many $\nu \in \mathfrak{X}_{\ell}$ we have
$$v_{p}\big{(} \cE(f,\eta\nu)\log_{\omega_{B_{f}}}(P_{f}(\eta \nu))\big{)}=0.$$
Here $\cE(f,\eta\nu)=1 - (\eta\nu)^{-1}(\overline{\frak{p}})\bfa_{p}(f)p^{-1}+(\eta\nu)^{-2}(\overline{\frak{p}})\epsilon(p)p^{-1}$.
\end{thm}
\begin{proof} 
In view of Corollary 4.2, it suffices to prove the non-triviality of the toric periods $P_{d^{-1}(f^{(p)}),\eta^{-1}}(\nu^{-1},n)$
modulo $p$.\\
\\ 
Being of weight zero, $d^{-1}(f^{(p)})$ is not a nearly holomorphic modular form (\cf \S2.2.3) 
and thus Theorem 3.1 does not seem to be directly applicable. The approach instead consists of two steps. Via an elementary congruence, we reduce the non-triviality to a setup involving the one for toric periods of a nearly holomorphic form. We then verify the non-triviality criteria in Theorem 3.1 for the nearly holomorphic form.\\
\\
\underline{Congruence}.
In view of the effect of the Katz $p$-adic differential operator on the $q$-expansion (\S2.2.3) and the $q$-expansion principle for $p$-adic modular forms (\S2.2.2), we have 
$$d^{(p-2)p}(f^{(p)})\equiv d^{-1}(f^{(p)}) \mod{\mathfrak{m}_p}$$ and   
$$d^{(p-2)p}(f^{(p)})\equiv d^{(p-2)p}(f) \mod{\mathfrak{m}_p}.$$ 
It thus follows that 
\begin{align*}
\sum_{[t]_{n} \in \Gamma_{n}^-} \phi([t]_{n}) d^{(p-2)p}(f)(x_{n}(t))({\bf{N}}_{K}^{(p-1)^2}\eta)(t)\\
\equiv \sum_{[t]_{n} \in \Gamma_{n}^-} \phi([t]_{n}) d^{-1}(f^{(p)})(x_{n}(t))\eta(t) \space \mod{\mathfrak{m}_p} 
\end{align*}
for $\phi: \Gamma_{n}^- \rightarrow \overline{\Z}_p$ as above. 
As $d^{(p-2)p}(f)$ is nearly holomorphic of weight $2(p-1)^{2}$, the pair $(d^{(p-2)p}(f),\eta {\bf{N}}_{K}^{(p-1)^2})$ satisfies the hypothesis (MC) 
and the toric periods are thus well defined.\\
\\
The above congruence can now be rewritten as 
\beq
P_{d^{(p-2)p}(f),{\bf{N}}_{K}^{(p-1)^2}\eta}(\phi,n)
\equiv P_ {d^{-1}(f^{(p)}),\eta}(\phi,n) \mod{\mathfrak{m}_p}.
\eeq
\noindent\\
\underline{Non-triviality}.
We may now apply Theorem 3.1 to the pair $(d^{(p-2)p}(f),\eta {\bf{N}}_{K}^{(p-1)^2})$. 
It thus suffices to verify the following.\\
\\
(H') For every $u \in \Z$ prime to $\ell$ and a given positive integer $r$, there exists a positive integer $\beta \equiv u \mod{\ell^{r}}$ and $a \in  {\bf{A}}_{K,f}^{\times}$ such that 
$v_{p}(\bfa_{\beta}(d^{(p-2)p}(f)^{\mathcal{R}},{\bf{c}}(a))) = 0$.\\
\\
Under the hypotheses (1)-(3), this is essentially verified in \cite[\S7.4]{Hs3}.
We give a brief summary in the appendix (\cf \S A.2).\\
\\
In view of the congruence (4.8), this finishes the proof.\\
\\
As non-triviality of the $p$-adic formal group logarithm of a point implies the point being non-torsion, ``In particular" part of the Theorem follows readily.\\
\end{proof}
\noindent\\
\begin{remark} 
The hypothesis $p$ being odd arises from the corresponding hypothesis in \cite{Hip}. Upon availability of Chai--Oort rigidity principle, it seems removable. 
The hypotheses (1)-(3) intervene only in the automorphic computation of the Fourier coefficients of the toric form.
The hypothesis (2) is perhaps removable. 
However, the hypotheses (1) and (3) are probably essential. We refer to the appendix (\cf Proposition A.4)
for a related discussion. 
\end{remark}
\noindent\\
\\
\subsection{Generalised Heegner cycles} 
\noindent In this section, we describe the basic setup regarding generalised Heegner cycles.\\
\\
\subsubsection{Generalities} 
In this subsection, we briefly recall generalities regarding generalised Heegner cycles following \cite[\S2]{BDP1} and \cite[\S4]{BDP3}.\\
\\
Unless otherwise stated, let the notations and hypotheses be as in the introduction. In particular, $p$ is an odd prime, $N$ a positive integer with $p \nmid N$, $b$ a positive integer with $(b, pN)=1$.
 Moreover, $K$ is an imaginary quadratic field with $p$ split satisfying the Heegner hypothesis (Hg) and $\cO$ the integer ring.
\\
\\
Recall that $X_1(N)$ is the modular curve of level $\Gamma_1(N)$ and $\infty$ the standard cusp $i\infty$. 
Here $i$ is a chosen square root of $-1$. 
Strictly speaking, the modular curve only exists as a Deligne--Mumford stack for $N \leq 3$. 
Let $r$ be a non-negative integer. 
Let $W_r$ be the $r$-fold Kuga--Sato variety over $X_{1}(N)$ constructed in \cite[App.]{BDP1}. 
In other words, $W_r$ is the canonical desingularisation of the $r$-fold self-product of the universal elliptic curve over $X_{1}(N)$.\\
\\
Let $A$ be a CM elliptic curve with CM by $\cO$ defined over the Hilbert class field $H$ (\cf \S2.1.4, \cite[\S1.4]{BDP1}). 
Recall $b$ denotes a positive integer prime to $N$ as above. 
For an ideal class $[\mathfrak{a}] \in \Pic(\cO_{bN})$ with $\mathfrak{a}$ prime to $N$, let $\varphi_{\mathfrak{a}}$ 
be the natural isogeny
\beq
\varphi_{\mathfrak{a}} : A \rightarrow A_{\mathfrak{a}}=A/A[\mathfrak{a}].
\eeq 
Strictly speaking, the above definition of the CM elliptic curve $A_{\mathfrak{a}}$ is correct only when $b=1$ and 
a suitable $\cO$-transform needs to be considered in the general case (\cf \cite[\S1.4]{BDP1}).\\
\\
Let $r_{1}$ and $r_2$ be non-negative integers such that $r_{1} \geq r_{2}$ and 
$r_{1} \equiv r_{2} \mod{2}$. Let $s$ and $u$ be non-negative integers such that
$$
r_{1}+r_{2}=2s \text{ and $r_{1}-r_{2}=2u.$}
$$
Let $X_{r_{1},r_{2}}$ be a $r_{1}+r_{2}+1$-dimensional variety given by 
$$X_{r_{1},r_{2}}=W_{r_{1}} \times A^{r_{2}}.$$ 
For an ideal $\mathfrak{a}$ as above, we consider 
$$ (A_{\mathfrak{a}} \times A)^{r_{2}} \times (A_{\mathfrak{a}} \times A_{\mathfrak{a}})^{u} = (A_{\mathfrak{a}}^{r_{1}} \times A^{r_{2}}) \subset    X_{r_{1},r_{2}}.          $$
The last inclusion arises from 
the embedding of $A_{\mathfrak{a}}^{r_{1}}$ in $W_{r_{1}}$ as a fiber of the natural projection $W_{r_{1}} \rightarrow X_{1}(N)$.\\
\\
Let 
$\Gamma_{\mathfrak{a}} \in Z^{1}(A_{\mathfrak{a}} \times A_{\mathfrak{a}})$ be the transpose of the
graph of $\sqrt{-d_{K}}$. 
Here $\sqrt{-d_{K}}$ is the square root of $-d_{K}$ whose image under the complex embedding $\iota_\infty$ 
has positive imaginary part. 
Let $\Gamma_{\varphi,\mathfrak{a}} \in Z^{1}(A_{\mathfrak{a}} \times A)$ be the transpose of the
graph of $\varphi_{\mathfrak{a}}$. 
Let 
$$   \Gamma_{r_{1},r_{2},\mathfrak{a}}  =     \Gamma_{\varphi,\mathfrak{a}}^{r_{2}} \times       \Gamma_{\mathfrak{a}}^{u}           $$
and 
\beq
\Delta_{r_{1},r_{2},\mathfrak{a}}=\epsilon_{X_{r_{1},r_{2}}}(\Gamma_{r_{1},r_{2},\mathfrak{a}}) \in CH^{s+1}(X_{r_{1},r_{2}}\otimes L)_{0,\Q}.
\eeq 
Here $\epsilon_{X_{r_{1},r_{2}}}$ is the idempotent in the ring of correspondences on $X_{r_{1},r_{2}}$ defined in \cite[\S4.1]{BDP3} 
and $L$ is the field of definition of the cycle $\Gamma_{r_{1},r_{2},\mathfrak{a}}$. 
The idempotent has an effect of making the cycle null-homologous (\cf \cite[Prop. 4.1.1]{BDP3}). 
The notation $CH^{s+1}(X_{r_{1},r_{2}}\otimes L)_{0,\Q}$ denotes the Chow group of homologically trivial cycles over $L$ of codimension $s+1$ with rational coefficients.\\
\\
When $r_{1}=r_{2}=r$, we let $X_{r}$, $\Gamma_{r,\mathfrak{a}}$ and $\Delta_{r,\mathfrak{a}}$ denote 
$X_{r,r}$, $\Gamma_{r, r,\mathfrak{a}}$ and $\Delta_{r, r,\mathfrak{a}}$, respectively.\\
\\
When $r_{1}$ is even and $r_{2}=0$, the above cycles are nothing but the classical Heegner cycles (\cf \cite[\S2.4]{BDP1}). 
Aspects of classical Heegner cycles were first investigated in \cite{Sc} and \cite{Ne}.\\
\\
When $r=0$, the generalised Heegner cycle $\Delta_{0,\mathfrak{a}}$ is a CM point on the modular curve $X_1(N)$. 
In this case, we replace $\Delta_{0,\mathfrak{a}}$ by $\Delta_{0,\mathfrak{a}} - \infty$ to make it homologically trivial. 
Note that this preserves the field of definition of the cycles as the $\infty$-cusp is defined over $\Q$.\\
\\
\\
\subsubsection{$p$-adic Abel--Jacobi map} 
In this subsection, we briefly recall generalities regarding a relevant $p$-adic Abel--Jacobi map 
following \cite[\S3]{BDP1}. 
For the dual exponential map, we refer to \cite{BK}.\\
\\
Let the notations and hypotheses be as in \S4.2.1. Let $\epsilon_X$ denote the idempotent $\epsilon_{X_{r_{1},r_{2}}}$.\\
\\
Let $F$ be a number field containing the Hilbert class field $H$. Let $V_s$ be the $p$-adic Galois representation of $G_F$ given by 
$$     V_{s}=H^{2s+1}_{\acute{e}t}(X_{r_{1},r_{2}} \times_{F} \overline{\Q}, \Q_p)                          .$$
Let $AJ_{F}^{\acute{e}t}$ be the étale Abel-Jacobi map
\beq
AJ_{F}^{\acute{e}t} : CH^{s+1}(X_{r_{1},r_{2}/F})_{0,\Q} \rightarrow H^{1}(F, \epsilon_{X} V_{s}(s+1))
\eeq
due to Bloch (\cf \cite{Ne} and \cite[Def. 3.1]{BDP1}). 
The Bloch--Kato conjecture implies that the $\Q_{p}$-linearisation $AJ_{F}^{\acute{e}t} \otimes_{\Q} \Q_p$ is injective (\cf \cite[(2.1)]{Ne}).\\
\\
Let $v$ be a place in $F$ above $p$ induced by the $p$-adic embedding $\iota_p$. We have the localisation map 
$$  
loc_{v}:     H^{1}(F, \epsilon_{X} V_{s}(s+1)) \rightarrow H^{1}(F_{v}, \epsilon_{X} V_{s}(s+1))
$$ 
in Galois cohomology.\\
\\
In general, the localisation map need not be injective.\\
\\
The image of the composition $loc_{v} \circ AJ_{F}^{\acute{e}t}$ 
 is contained in the Bloch--Kato subgroup $H^{1}_{f}(F_{v}, \epsilon_{X} V_{s}(s+1))$ (\cf \cite[Thm. 3.1 (ii)]{Ne} and \cite{Sa}). 
In terms of the interpretation of the local cohomology $H^{1}(F_{v}, \epsilon_{X} V_{s}(s+1))$
as a group of extension classes, the elements of the subgroup correspond to as crystalline extensions 
as $\epsilon_{X} V_{s}(s+1)$ lets them to be.\\
\\
In view of the Bloch--Kato logarithm map and de Rham Poincaré duality, 
we have a canonical isomorphism
\beq
   \log :  H^{1}_{f}(F_{v}, \epsilon_{X} V_{s}(s+1)) \iso (Fil^{s+1}\epsilon_{X}H^{2s+1}_{dR}(X_{r_{1},r_{2}/F_{v}})(s))^{\vee}                           
\eeq
(\cf \cite[\S3.4]{BDP1}). 
Here $Fil^{\bullet}$ is the Hodge filtration on $\epsilon_{X}H^{2r+1}_{dR}(X_{r/F_{v}})(r)$ and 
$(Fil^{r+1}\epsilon_{X}H^{2r+1}_{dR}(X_{r/F_{v}})(r))^{\vee}$ the dual arising 
from the Poincaré pairing.\\
\\
The composition with the étale Abel--Jacobi map gives rise to the $p$-adic Abel--Jacobi map
\beq
AJ_{F}: CH^{s+1}(X_{r_{1},r_{2}/F})_{0,\Q} \rightarrow (Fil^{s+1}\epsilon_{X}H^{2s+1}_{dR}(X_{r_{1},r_{2}/F_{v}})(s))^{\vee}. 
\eeq 
There does not seem to be a general conjecture regarding the image and kernel of the $p$-adic Abel--Jacobi map. 
However,\\
\\
(BBK) the Bloch--Beilinson and the Bloch--Kato conjectures suggest to investigate contexts where the $\Q_p$-linearisation $AJ_{F} \otimes_{\Q} \Q_p$ is an isomorphism 
or generically so (\cf \cite[(2.1)]{Ne}).\\
\\
In general, the $\Q_p$-linearisation $AJ_{F} \otimes_{\Q} \Q_p$ need not be injective or surjective. 
For a relation of the non-triviality of the $p$-adic Abel--Jacobi map 
to coniveau filtration and refined Bloch--Beilinson conjecture, we refer to \cite[\S2]{BDP3}.\\
\\
For later purposes, we recall an explicit description of the middle step  
$Fil^{s+1}\epsilon_{X}H^{2s+1}_{dR}(X_{r_{1},r_{2}/F_{v}})(s)$ of the Hodge filtration.\\
\\
For a positive integer $k\geq 2$, let 
$S_{k}(\Gamma_{1}(N), F_v)$ denote the space of cusp form of level $\Gamma_{1}(N)$ with coefficients in $F_v$. 
We have a canonical isomorphism 
\beq
S_{r_{1}+2}(\Gamma_{1}(N),F_{v}) \otimes_{F_{v}} \Sym^{r_{2}} H^{1}_{dR}(A_{/F_{v}}) \iso Fil^{s+1}\epsilon_{X}H^{2s+1}_{dR}(X_{r_{1},r_{2}/F_{v}})(s),
\eeq 
essentially due to Scholl and it arises from 
$$  f \otimes \alpha \mapsto \omega_{f} \wedge \alpha.                             $$ 
Here $\omega_f$ is the normalised differential associated with $f \in S_{r_{1}+2}(\Gamma_{1}(N),F_{v})$ (for example,  
\cite[Cor 2.3]{BDP1}) and $\alpha \in \Sym^{r_{2}} H^{1}_{dR}(A/F_{v})$. 
We refer to \cite[\S2.2]{BDP1} for details.\\
\\
The symmetric power $\Sym^{r_{2}} H^{1}_{dR}(A/F_{v})$ in turn has the following basis. 
Let $\omega_{A} \in \Omega^{1}_{A/F_{v}}$ be a non-zero differential. Let $\eta_{A} \in \Omega^{1}_{A/F_{v}}$ 
be the corresponding element in \cite[(1.4.2)]{BDP1}. It may be directly seen that 
$\{\omega_{A}, \eta_{A}\}$ is a basis of $\Omega^{1}_{A/F_{v}}$. 
Thus, $\big{\{} \omega_{A}^{j} \eta_{A}^{r_{2}-j}: 0 \leq j \leq r_{2} \big{\}}$ is a basis of $\Sym^{r_{2}} H^{1}_{dR}(A/F_{v})$.\\
\\
For $0 \leq j \leq r_{2}$ and $f$ varying over Hecke eigenforms in $S_{r_{1}+2}(\Gamma_{1}(N),F_{v})$, 
the wedge products $$\bigg{\{}\omega_{f} \wedge \omega_{A}^{j} \eta_{A}^{r_{2}-j}: f , 0 \leq j \leq r_{2} \bigg{\}}$$ thus form a basis 
of $Fil^{s+1}\epsilon_{X}H^{2s+1}_{dR}(X_{r_{1},r_{2}/F_{v}})(s)$.\\
\\
\\
\subsubsection{Conjectures} In this subsection, we state conjectures
regarding the non-triviality of generalised Heegner cycles and the image under 
the étale and $p$-adic Abel--Jacobi maps over anticyclotomic extensions 
of an imaginary quadratic extension. 
These conjectures are essentially an ammendment of the conjectures due to Mazur, Bloch--Beilinson and 
Bloch--Kato. 
For the latter conjectures, we refer to \cite{M} and \cite{BK}.\\
\\
Let the notations and hypotheses be as in \S4.2.1. 
In particular, $N$ is a positive integer such that $p\ndivide N$. 
For a positive integer $k\geq 2$, let $S_{k}(\Gamma_{0}(N),\epsilon)$ be the space of elliptic cusp forms of weight $k$, level $\Gamma_{0}(N)$ and neben-character $\epsilon$.
Let $f\in S_{r_{1}+2}(\Gamma_{0}(N),\epsilon)$ be a normalised newform. 
In particular, it is a Hecke eigenform with respect to all Hecke operators. 
Let $N_{\epsilon}|N$ be the conductor of $\epsilon$. 
Let $\mathfrak{N}_{\epsilon}|\mathfrak{N}$ be the unique ideal of norm $N_{\epsilon}$. 
The existence follows from the generalised Heegner hypothesis (Hg).\\
\\
Let ${\bf{N}}: {\bf{A}}_{\Q}^{\times}/\Q^{\times} \rightarrow \C^\times$ be the norm Hecke character over $\Q$ given by 
$$     
{\bf{N}}(x)=||x||.        
$$
Here $||\cdot||$ denotes the adelic norm. 
Let ${\bf{N}}_{K}:={\bf{N}}\circ N_{\Q}^{K}$ be the norm Hecke character over $K$ for the relative norm $N_{\Q}^{K}$.
For a Hecke character $\lam:{\bf{A}}_{K}^{\times}/K^{\times} \rightarrow \C^\times$ over $K$, let $\mathfrak{f}_\lam$ (resp. $\epsilon_\lam$) denote 
its conductor (resp. the restriction $\lam|_{{\bf{A}}_{\Q}^\times}$) as before.\\
\\
Let $b$ be a positive integer prime to $N$. 
Let $\Sg_{r_{1},r_{2}}(b,\mathfrak{N},\epsilon)$ be the set of Hecke characters $\lam$ such that:\\
\\
(C1) $\lam$ is of infinity type $(j_{1},j_{2})$ with $j_{1}+j_{2}=r_{2}$ and $\epsilon_{\lam}=\epsilon {\bf{N}}^{r_{2}}$,\\
(C2) $\mathfrak{f}_{\lam}=b \cdot \mathfrak{N}_\epsilon$ and\\
(C3) we have $\epsilon_{q}(f,\lam^{-1}{\bf{N}}_{K}^{-u})=1$ for $\epsilon_{q}(f,\lam^{-1}{\bf{N}}_{K}^{-u})$ the local root number corresponding to Rankin--Selberg convolution of the pair $(f,\lam^{-1}{\bf{N}}_{K}^{-u})$ for all finite primes $q$.\\
\\
Let $\Sg_{r_{1},r_{2}}^{(2)}(b,\mathfrak{N},\epsilon)$ be the subset of $\Sg_{r_{1},r_{2}}(b,\mathfrak{N},\epsilon)$ 
whose elements have infinity type $(r_{2}+1-j, 1+j)$ for some $0 \leq j \leq r_2$.\\
\\
The underlying global root numers are given by the following
\begin{lm}\label{RN} 
Let $f\in S_{r_{1}+2}(\Gamma_{0}(N),\epsilon)$ be a normalised newform for an integer $r_{1} \geq 0$. Let $K$ be an imaginary quadratic field satisfying Heegner hypothesis (Hg). For a positive integer $b$ prime to $N$, let $\Sg_{r_{1},r_{2}}(b,\mathfrak{N},\epsilon)$ be a set of arithmetic Hecke characters over $K$ satisfying (C1)-(C3) and $\Sg_{r_{1},r_{2}}^{(2)}(b,\mathfrak{N},\epsilon)$ its subset as above. Let $\chi \in \Sg_{r_{1},r_{2}}(b,\mathfrak{N},\epsilon)$, then 
$$
\epsilon(f,\chi^{-1}{\bf{N}}_{K}^{-u})=-1 \text{(resp. $1$)}
$$
for $\chi \in \Sg_{r_{1},r_{2}}^{(2)}(b,\mathfrak{N},\epsilon)$ (resp. else). 
Here $\epsilon(f,\chi^{-1}{\bf{N}}_{K}^{-u})$ denotes the global root number corresponding to Rankin--Selberg convolution of the pair $(f,\chi^{-1}{\bf{N}}_{K}^{-u})$ (\cf \cite[\S4.1]{BDP1}).
\end{lm}
\noindent\\
Let $\chi \in \Sg_{r_{1},r_{2}}^{(2)}(b,\mathfrak{N},\epsilon)$ with infinity type $(r_{2}+1-j, 1+j)$ for some $0 \leq j \leq r_2$.\\
\\
Let $M_f$ (resp. $M_{\chi^{-1}}$) be the Grothendieck motive associated with $f$ (resp. $\chi^{-1}$). 
By the K$\ddot{u}$nneth formula, the motive $H^{s+1}(X_{r_{1},r_{2}})\otimes_{\Q} L$ contains 
$M_{f} \otimes_{\Q} M_{\chi^{-1}}$ as a submotive 
for a sufficiently large number field $L$. 
For a minimal choice of $L$, we refer to \cite[\S4.2]{BDP3}.
Let $\pi_{f,\chi}$ be the projector defining the submotive (\cf \cite{Sc}).\\
\\
Let $S_b$ be a set of representatives for $\Pic(\cO_{bN})$ consisting of ideals prime to $N$. 
We now regard $\chi$ as an ideal Hecke character. 
Let $H_{\chi}$ be the abelian extension of $H$ generated by the values of $\chi$ on $S_b$. 
In particular, the extension $H_\chi$ depends on the choice of $S_b$. 
Let $v=v_{\chi}$ be the place above $p$ in $H_{\chi}$ induced via the $p$-adic embedding $\iota_p$.\\
\\
Let 
\beq
 \Delta_{\chi} = \sum_{[\mathfrak{a}] \in S_{b}} \chi^{-1}(\mathfrak{a}) {\bf{N}}(\mathfrak{a}) \Delta_{\varphi_{\mathfrak{a}}}  \in CH^{s+1}(X_{r_{1},r_{2}/H_{\chi}})_{0,\Q} \otimes L.                                                             
\eeq 
The cycle depends on the choice of representatives $S_b$. Moreover, 
it is defined over the extension $H_{\chi}$ by the definition.\\
\begin{defn} 
The generalised Heegner cycle $\Delta_{f,\chi}$ associated with the pair $(f,\chi)$ is given by
$$ \Delta_{f,\chi}= \pi_{f,\chi}(\Delta_{\chi}).         $$
\end{defn}
\noindent The generalised Heegner cycle is independent of the choice of representatives $S_b$
up to a cycle in the kernel of the complex Abel--Jacobi map (\cf \cite[Rem. 4.2.4]{BDP3}).\\
\\
We consider the non-triviality of the cycles $\Delta_{f,\chi}$, as $\chi$ varies.\\
\\
Under the hypotheses (C1)-(C3), 
the Rankin--Selberg $L$-function 
corresponding to the pair $(f, \chi^{-1}{\bf{N}}_{K}^{-u})$ is self-dual with root number $-1$ (\cf Lemma \ref{RN}). 
The Bloch--Beilinson conjecture accordingly implies the existence of a non-torsion 
homologically trivial cycle in the Chow realisation of the motive 
$M_{f} \otimes_{\Q} M_{\chi^{-1}{\bf{N}}_{K}^{-u}}$. 
Perhaps, a natural candidate is the cycle $\Delta_{f,\chi{\bf{N}}_{K}^{u}}$. 
There exists an algebraic correspondence from $X_{r_{1}}$ to $X_{r_{1},r_{2}}$ mapping the cycle $\Delta_{f,\chi{\bf{N}}_{K}^{u}}$ to the $\Delta_{f,\chi}$ (\cf \cite[Prop. 4.1.1]{BDP3}). 
One may thus expect a generic non-triviality of the cycles $\Delta_{f,\chi}$, as $\chi$ varies.\\
\\
Based on Mazur's consideration in the case of weight two (\cf \cite{M}), an Iwasawa theoretic family of the cycles arises as follows. 
We first fix a Hecke character $\eta \in \Sg_{r_{1},r_{2}}^{(2)}(c,\mathfrak{N},\epsilon)$, for some $c$ prime to $N$. 
Let $\ell$ be an odd prime unramified in $K$ and prime to $cN$. 
Let $H_{cN\ell^\infty} = \bigcup_{n\geq 0} H_{cN\ell^n}$ be the ring class field of conductor $cN\ell^\infty$. 
Let $K_{\infty}^- \subset H_{cN\ell^\infty}$ be the anticyclotomic $\Z_\ell$-extension of $K$.
Let $G_{n}=\Gal(H_{cN\ell^n}/K)$
and $\Gamma_{\ell}=\varprojlim G_{n}$. 
Let $\mathfrak{X}_{\ell}$ be the subgroup of all characters of finite order of the group $\Gal(K_{\infty}^-/K) \iso \Z_\ell$. 
As $\nu \in \mathfrak{X}_{\ell}$ varies, we consider non-triviality of the generalised Heegner cycles $\Delta_{f,\eta\nu}$.\\
\\
\\
{\bf{Conjecture A.}} 
Let $f\in S_{r_{1}+2}(\Gamma_{0}(N),\epsilon)$ be a normalised newform for an integer $r_{1} \geq 0$. Let $K$ be an imaginary quadratic field satisfying Heegner hypothesis (Hg). For a positive integer $b$ prime to $N$, let $\Sg_{r_{1},r_{2}}^{(2)}(b,\mathfrak{N},\epsilon)$ be a set of arithmetic Hecke characters over $K$ satisfying (C1)-(C3) as above. 
Let $\eta \in \Sg_{r_{1},r_{2}}^{(2)}(c,\mathfrak{N},\epsilon)$, for some $c$ prime to $N$. 
Let $\ell$ be an odd prime unramified in $K$ and prime to $cN$ and $\mathfrak{X}_\ell$ the set of finite order Hecke characters with $\ell$-power order as above.
Then, for all but finitely many $\nu \in \mathfrak{X}_{\ell}$ we have 
$$     \Delta_{f,\eta\nu} \neq 0                $$
for $\Delta_{f,\eta\nu}$ a generalised Heegner cycle as above.\\
\noindent\\
\\
We also consider the non-triviality of the generalised Heegner cycles under the étale Abel--Jacobi map.\\
\\
Let $p$ be an odd prime relatively prime to $Nc$. 
Let $M_{f,\acute{e}t}$ (resp. $M_{\chi^{-1},\acute{e}t}$) be the $p$-adic étale realisation of the motive $M_{f}$ (resp. $M_{\chi^{-1}}$). 
We accordingly have
$$  AJ_{H_{\eta\nu},\acute{e}t}(\Delta_{f,\chi}) \in H^{1}(H_{\eta\nu},\epsilon_{X}( M_{f,\acute{e}t} \otimes_{\Q}  M_{\chi^{-1}})).                $$
\\
{\bf{Conjecture B.}} 
Let $f\in S_{r_{1}+2}(\Gamma_{0}(N),\epsilon)$ be a normalised newform for an integer $r_{1} \geq 0$. Let $K$ be an imaginary quadratic field satisfying Heegner hypothesis (Hg). For a positive integer $b$ prime to $N$, let $\Sg_{r_{1},r_{2}}^{(2)}(b,\mathfrak{N},\epsilon)$ be a set of arithmetic Hecke characters over $K$ satisfying (C1)-(C3) as above. 
Let $\eta \in \Sg_{r_{1},r_{2}}^{(2)}(c,\mathfrak{N},\epsilon)$, for some $c$ prime to $N$. 
Let $\ell$ be an odd prime unramified in $K$ and prime to $cN$ and $\mathfrak{X}_\ell$ the set of finite order Hecke characters with $\ell$-power order as above. 
Let $p$ be an odd prime relatively prime to $Nc$.
Then, for all but finitely many $\nu \in \mathfrak{X}_{\ell}$ we have 
$$     AJ_{H_{\eta\nu},\acute{e}t}(\Delta_{f,\eta\nu}) \neq 0                $$
for $\Delta_{f,\eta\nu}$ a generalised Heegner cycle as above and $AJ_{H_{\eta\nu},\acute{e}t}(\cdot)$ the \'etale Abel--Jacobi image.\\
\noindent\\
We finally consider the non-triviality of the generalised Heegner cycles under the $p$-adic Abel--Jacobi map.\\
\\
We first recall that the Bloch--Kato subgroup $H_{f}^{1}(H_{\eta\nu,v},\epsilon_{X}( M_{f,\acute{e}t} \otimes_{\Q}  M_{\chi^{-1}}))$ is one-dimensional 
over $H_{\eta\nu,v}$. 
Moreover, under the Bloch--Kato logarithm map 
the corresponding basis is given by $\omega_{f} \wedge \omega_{A}^{j}\eta_{A}^{r_{2}-j}$ (\cf \S4.2.2).\\
\\
\\
{\bf{Conjecture C.}} 
Let $f\in S_{r_{1}+2}(\Gamma_{0}(N),\epsilon)$ be a normalised newform for an integer $r_{1} \geq 0$. Let $K$ be an imaginary quadratic field satisfying Heegner hypothesis (Hg). For a positive integer $b$ prime to $N$, let $\Sg_{r_{1},r_{2}}^{(2)}(b,\mathfrak{N},\epsilon)$ be a set of arithmetic Hecke characters over $K$ satisfying (C1)-(C3) as above. 
Let $\eta \in \Sg_{r_{1},r_{2}}^{(2)}(c,\mathfrak{N},\epsilon)$, for some $c$ prime to $N$. 
Let $\ell$ be an odd prime unramified in $K$ and prime to $cN$ and $\mathfrak{X}_\ell$ the set of finite order Hecke characters with $\ell$-power order as above. 
Let $p$ be an odd prime relatively prime to $Nc$.
Then, for all but finitely many $\nu \in \mathfrak{X}_{\ell}$ we have 
$$     AJ_{H_{\eta\nu}}(\Delta_{f,\eta\nu})(\omega_{f} \wedge \omega_{A}^{j}\eta_{A}^{r_{2}-j}) \neq 0                $$
for $\Delta_{f,\eta\nu}$ a generalised Heegner cycle as above and $AJ_{H_{\eta\nu},\acute{e}t}(\cdot)(\omega_{f} \wedge \omega_{A}^{j}\eta_{A}^{r_{2}-j})$ the $p$-adic Abel--Jacobi image.
\noindent\\
\\
Conjecture C (resp. Conjecture B) evidently implies Conjecture B (resp. Conjecture A).\\
\\
We remind that Conjecture C is in accordance with the Bloch--Beilinson and the Bloch--Kato conjecture (\cf (BBK)).\\
\\
As a refinement of Conjecture C, we can ask for a generic mod $p$ non-vanishing of normalised $p$-adic Abel--Jacobi image. 
We do not precisely formulate the refinement. 
In this article, we instead describe results towards the refinement in the case $l \neq p$ and under hypothesis (ord).\\
\\
\\
\subsection{Non-triviality, II}
In this subsection, we prove the non-triviality of the $p$-adic Abel--Jacobi image of generalised Heegner cycles modulo $p$ (\cf Conjecture C). 
As the details of our approach are quite similar to the ones in \S4.1, 
we outline a sketch.\\
\\
Unless otherwise stated, let the notations and hypotheses be as in \S4.2. 
In particular, $p$ is an odd prime, $N$ a positive integer with $p \nmid N$, $b$ a positive integer with $(b, pN)=1$ and $\ell$ a prime such that $\ell \nmid 2pNb$. Moreover, $K$ is an imaginary quadratic field with $p$ split satisfying the Heegner hypothesis (Hg) with $\ell$ unramified, 
$K_{\ell^n}^{-}$ is the anticyclotomic extension of $K$ with conductor $bN \cdot \ell^n$, 
$\Gamma_{n}^{-}=\Gal(K_{\ell^n}^{-}/K)$ the anticyclotomic Galois group and $\Gamma_{\ell}^{-}=\varprojlim_{n} \Gamma_{n}^-$ and $\mathfrak{X}_\ell$ the set of finite order characters of $\Gamma_\ell$.\\
\\
Further, $f$ is a normalised newform of weight at least $2$, level $\Gamma_{0}(N)$ and $\eta$ a Hecke character over $K$ with conductor divisible by $b$ and weight of $f$ being dominant such that the Rankin--Selberg convolution corresponding to the pair $(f,\eta)$ is self-dual with root number $-1$. 
For $\nu \in \mathfrak{X}_\ell$, the Rankin--Selberg convolution corresponding to the pair $(f,\eta\nu)$ is again self-dual with root number $-1$ (Lemma 4.4). 
We consider non-triviality of the $p$-adic Abel--Jacobi image of generalised Heegner cycles $\Delta_{f,\eta\nu}$ as $\nu \in \mathfrak{X}_{\ell}$ varies.\\
\\
In the residually non-CM case, our result towards a refinement of Conjecture C is the following.\\
\begin{thm}
Let $p$ be an odd prime and $N$ a positive integer such that $p \ndivide N$. 
Let $K$ be an imaginary quadratic field satisfying (ord) and (Hg) with $\mathfrak{p}$ a prime above $p$ determined via an initial embedding $\iota_p: \overline{\Q} \hookrightarrow \C_p$. 
Let $f\in S_{r_{1}+2}(\Gamma_{0}(N),\epsilon)$ be a normalised newform with $T_{p}$-eigenvalue $\bfa_{p}(f)$ and $\eta \in \Sg_{r_{1},r_{2}}^{(2)}(c,\mathfrak{N},\epsilon)$ a Hecke character of infinity type $(r_{2}+1-j,1+j)$, 
for some $c$ prime to $N$ and $0 \leq j \leq r_2$. 
Suppose that\\
\\
(1). The residual representation $\rho_{f}|_{G_{K}} \mod{\mathfrak{m}_{p}}$ is absolutely irreducible,\\
(2). $N^-$ is square-free and\\
(3). $p \ndivide \prod_{v|c^{-}}|\Delta_{\eta,v}|$ for $\Delta_{\eta,v}$ the finite group $\eta(\cO_{K_{v}}^{\times})$.\\
\\
Let $\ell\neq p$ be an odd prime unramified in $K$ and prime to $cN$. 
Let $\mathfrak{X}_\ell$ be the set of $\ell$-power order anticyclotomic Hecke characters over $K$ as above. 
Then, for all but finitely many $\nu \in \mathfrak{X}_{\ell}$ we have
$$v_{p}\bigg{(}\frac{\cE_{p}(f,\eta\nu{\bf{N}}_{K}^{u})}{j!} AJ_{H_{\eta\nu}}(\Delta_{f,\eta\nu})(\omega_{f} \wedge \omega_{A}^{j}\eta_{A}^{r-j})\bigg{)}=0.$$
\noindent\\
Here $\cE_{p}(f,\eta\nu{\bf{N}}_{K}^{u})=1-(\eta\nu{\bf{N}}_{K}^{u})^{-1}(\overline{\mathfrak{p}})\bfa_{p}(f)+(\eta\nu{\bf{N}}_{K}^{u})^{-2}(\overline{\mathfrak{p}})\epsilon(p)p^{r_{2}+1}$. 
\noindent\\
\end{thm}
\begin{proof} 
In view of \cite[Prop. 4.1.2]{BDP3}, it suffices to restrict to the case $r_{1}=r_{2}=r$.\\
\\
Let $f^{(p)}$ be the $p$-depletion of the $f$ as in (4.1). The pair $(d^{-1-j}(f^{(p)}),\eta {\bf{N}}_{K}^{-1-j})$ satisfies the hypothesis (MC) 
and the toric periods are well defined.\\
\\
\underline{$p$-adic Waldspurger.}
For $\nu$ factoring through the anticyclotomic quotient $\Gamma_n^-$, the toric period corresponding to the pair $(d^{-1-j}(f^{(p)}),\eta^{-1}\nu {\bf{N}}_{K}^{-1-j})$
 is given by  
\beq
P_{d^{-1-j}(f^{(p)}),\eta {\bf{N}}_{K}^{-1-j}}(\nu,n) = 
\frac{\cE_{p}(f,\eta\nu)}{j!} \sum_{[\mathfrak{a}]\in \Gamma_{n}} (\eta\nu)^{-1}(\mathfrak{a}){\bf{N}}_{K}(\mathfrak{a})\cdot AJ_{F}(\Delta_{\varphi_{\mathfrak{a}}\varphi_{0}})(\omega_{f} \wedge \omega_{A}^{j}\eta_{A}^{r-j}).
\eeq
This is again $p$-adic Waldspurger formula due to Bertolini--Darmon--Prasanna (\cite[\S3]{BDP1}) and the proof is similar to that of Theorem 4.1.\\
\\
For the non-triviality of these toric periods modulo $p$, we again first find a pair whose toric periods are congruent modulo $p$ to the periods of interest and
such that Theorem 3.1 can be applied to the new pair.\\
\\
\underline{Congruence.}
We consider a pair $(d^{m(p-2)p-j}(f),\eta {\bf{N}}_{K}^{m(p-2)p-j})$ for $m$ is a positive integer such that 
$$
p-1 \mid m(p-2)p+1, k+2m(p-2)p-2j \geq 2.
$$
As in the proof of Theorem 4.3, we have the congruence\\
\beq
P_{d^{m(p-2)p-j}(f),\eta {\bf{N}}_{K}^{m(p-2)p-j}}(\phi,n)
\equiv  P_{d^{-1-j}(f^{(p)}),\eta {\bf{N}}_{K}^{-1-j}}(\phi,n)\mod{\mathfrak{m}_p}.
\eeq
\noindent\\
The non-triviality hypothesis (H) can be verified for the pair in an analogous manner as in the proof of Theorem 4.3. We refer to the appendix (\cf Proposition A.4) for details\\
\end{proof}
\noindent\\
\begin{remark} 
(1). In view of \cite[\S2.3]{BDP3}, 
it follows that the generalised Heegner cycles are non-trivial in the top graded piece of the coniveau filtration on the 
middle-codimensional Chow group 
of $X_{r_{1},r_{2}}$ over the $\Z_l$-anticyclotomic extension of $K$. 
The non-triviality is thus an evidence for the refined Bloch--Beilinson conjecture (\cf \cite[\S2]{BDP3}).\\
\\
(2). In \cite{Ho}, Howard proved a weak version of analogous non-triviality of the étale Abel--Jacobi image of 
classical Heegner cycles for infinitely many $\nu \in \mathfrak{X}_\ell$ (\cf Conjecture B). 
Howard's approach is a cohomological adoption of Cornut's approach (\cf \cite{C}). The approach does not seem relevant for the $p$-adic Abel--Jacobi image as above.\\
\\
(3). It seems likely that the above approach would also work for the case when $f$ has CM over $K$.\\
\end{remark}
\noindent\\
\section{$p$-indivisibility of Heegner points}
\noindent In this section, we prove the $p$-indivisibility of the Heegner points over the $\Z_l$-anticyclotomic extension of $K$ when the underlying weight two Hecke-eigen cuspform is $p$-ordinary.
This is based on Theorem 4.3.\\
\\
Let the notation and hypotheses be as in the introduction. 
In particular, $p$ is an odd prime, $N$ a positive integer with $p \nmid N$, $b$ a positive integer with $(b, pN)=1$ and $\ell$ a prime such that $\ell \nmid 2pNb$. Moreover, $K$ is an imaginary quadratic field with $p$ split satisfying the Heegner hypothesis (Hg), 
$K_{\ell^n}^{-}$ is the anticyclotomic extension of $K$ with conductor $bN \cdot \ell^n$, 
$\Gamma_{n}^{-}=\Gal(K_{\ell^n}^{-}/K)$ the anticyclotomic Galois group and $\Gamma_{\ell}^{-}=\varprojlim_{n} \Gamma_{n}^-$.
Further, $\mathfrak{p}$ denotes a prime above $p$ in $K$ induced by the $p$-adic embedding $\iota_p$.
Let $\mathfrak{p}_n$ be a prime above $\mathfrak{p}$ in $K_{bNl^n}$  induced by $\iota_p$. Let $K_{n}/\Q_p$ denote 
the $p$-adic local field $K_{bN\ell^n,\mathfrak{p}_n}$ and $O_n$ the integer ring. 
By definition, $K_{n}/\Q_p$ is an unramified extension.\\
\\
Recall $f \in S_{2}(\Gamma_{0}(N),\epsilon)$ is a weight two Hecke-eigen cuspform with $T_{p}$-eigenvalue $\bfa_{p}(f)$ and $B_{f}$ an abelian variety associated to $f$. 
We suppose that $f$ is $p$-ordinary.\\
%
\\
Our result regarding the $p$-indivisibility of the Heegner points is as follows.\\
\begin{thm} 
Let $p$ be an odd prime and $N$ a positive integer such that $p \ndivide N$. 
Let $K$ be an imaginary quadratic field satisfying (ord) and (Hg) with $\mathfrak{p}$ a prime above $p$ determined via an initial embedding $\iota_p: \overline{\Q} \hookrightarrow \C_p$. 
Let $f\in S_{2}(\Gamma_{0}(N),\epsilon)$ be a $p$-ordinary normalised newform with $T_{p}$-eigenvalue $\bfa_{p}(f)$, $B_{f}$ a corresponding $p$-optimal abelian variety and $c$ a positive integer prime to $pN$.
Suppose that\\
\\
(1). The residual representation $\rho_{f}|_{G_{K}} \mod{\mathfrak{m}_{p}}$ is absolutely irreducible and\\
(2). $N^-$ is square-free.\\
\\
Let $\ell\neq p$ be an odd prime unramified in $K$ and prime to $cN$. 
Then, there exists a positive integer $n_0$ such that for all integers $n \geq n_0$, the Heegner points $\Phi_{f}(\Delta_{c\ell^n})$ 
are $p$-indivisible in the Mordell--Weil group $B_{f}(K_{n})$.\\
\end{thm}
\begin{proof} 
As $f$ is $p$-ordinary Hecke-eigen cuspform, the corresponding abelian variety $B_{f}$ is $p$-ordinary and extends as an abelian scheme over $O_n$.
It follows that 
for a non-torsion $P\in B_{f}(K_{n})$, we have $$v_{p}(\log_{\omega_{B_{f}}}(P))\geq 1.$$
\noindent
In view of Theorem 4.3, there exists an integer $n_0$ such that for all integers $n\geq n_0$, the Heegner points $\Phi_{f}(\Delta_{c\ell^n})$ are non-torsion. 
This concludes the proof.\\
\end{proof}
\noindent
\begin{remark}
(1). In \cite{GP}, arithmetic applications of the $p$-indivisibility of the Heegner point defined over the imaginary quadratic extension $K$ are considered under different hypothesis including the prime $\ell$ being inert in $K$. 
It seems likely that Theorem 5.1 also has similar arithmetic applications.\\
\\
(2). A version of Theorem 5.1 holds for a pair $(f,\chi)$ as in Theorem 4.3. We restrict to the above version for simplicity.\\
\end{remark}
\noindent\\
\appendix
\section{Toric modular forms}
\noindent In this appendix, we describe toric modular forms closely following the exposition in \cite{Hs3}. Toric forms constitute a key automorphic ingredient in the non-triviality of toric periods (\cf Hypothesis (H)).\\
\\
In the appendix, we mainly follow automorphic normalisations.\\

\subsection{Setup} In this subsection, we describe the basic setup.\\
\\
We begin with a few underlying notions. 
For a non-archimedean local field $F$, let $\mu,\nu:F^\x\to\C^\x$ be two characters of $F^\x$. Let $I(\mu,\nu)$ be the space consisting of smooth and $\GL_2(\cO_F)$-finite functions $f:\GL_2(F)\to\C$ such that
\[f(\MX{a}{b}{0}{d}g)=\mu(a)\nu(d)\abs{\frac{a}{d}}^\onehalf f(g).\]
$I(\mu,\nu)$ turns out to be an admissible representation of $\GL_2(F)$. Let $\pi(\mu,\nu)$ be the unique infinite dimensional subquotient of $I(\mu,\nu)$.  We call $\pi(\mu,\nu)$ a principal series if $\mu\nu^{-1}\not=\Abs^\pm$ and a special representation if $\mu\nu^{-1}=\Abs^\pm$.\\
\\
For a local field $F$, let $\pi$ be an irreducible admissible representation of $\GL_2(F)$ and $\addchar:F\to\C^\x$ a non-trivial additive character. Let $\cW(\pi,\addchar)$ be the Whittaker model of $\pi$ given by the subspace of smooth functions $W:\GL_2(F)\to\C$ such that
\begin{mylist}
\item $W(\MX{1}{x}{0}{1}g)=\addchar(x)W(g)$ for all $x\in F$.
\item If $v$ is archimedean, $W(\DII{a}{1})=O(\abs{a}^M)$ for some positive number $M$.
\end{mylist}
(\cf\cite[Thm.\,6.3]{JL}).\\
\\
We now turn towards the setup. 
Let the notation and hypothesis be as in \S4.2.3.
In particular, $(f,\chi)$ is a pair of an elliptic newform with conductor $N$ neben-type $\epsilon$ and a Hecke character over imaginary quadratic field $K$ Let $\pi$ be the cuspidal automorphic representation of $\GL_2(\A)$ generated by the Hecke eigenform $f$. For a place $v$ of $\Q$, let $\pi_v$ be the corresponding representation of $\GL_{2}(\Q_v)$. For $v \nmid N$, let $W_v^\circ$ be the new Whittaker function (\cite[\S3.6]{Hs3}).
\\
\\
As in the proof of Theorem 4.3, we suppose $r_{1}=r_{2}=r$ for $r \geq 0$. In particular, $f$ is with weight $r+2$ and $\chi$
with infinity type $(r+1-j,1+j)$ for some $0\leq j \leq r$.
\\
\\
Let $\frakC_\chi$ (resp. $N_\epsilon$) be the conductor of $\chi$ (resp. $\epsilon$) and $c_\chi=\frakC_\chi\cap\Q$. 
For the conductor $N$ of the newform $f$, we have $N=N^{+}N^-$ with $N^+$ precisely divisible by split primes in the extension $K/\Q$ as before.
We further decompose $N^-=N^-_s N^-_r$ with $N^-_r$ precisely divisible by prime factors of $N_\epsilon$.\\
\\
As before, let $\bdh$ denote the set of finite places of $\Q$. 
For $v \in \bdh$, let \beq\label{E:alv.W}\begin{aligned}c_v(\chi)=&\inf\stt{n\in\Z_{\geq 0}\mid \chi=1\text{ on }(1+\uf^n\cO_K)^\x},\\
m_v(\chi,\pi)=&c_v(\chi)-v(N^-).\end{aligned}\eeq 
Here $\uf$ denotes a uniformiser at $v$. 
We note that $c_v(\chi)=v(c_\chi)$. Let
\begin{align*}
A(\chi)=&\stt{v\in\bdh\mid\text{$K_v$ is a field, $\pi_v$ is special and $c_v(\chi)=0$}}
\end{align*}
and
\[B(\chi)=\stt{v\in\bdh\mid\text{$v$ is non-split in $K$ with $c_v(\chi)>0$}}.\]\\
For $v\in B(\chi)$, 
 say $\pi_{v}=\pi(\mu,\nu)$ with unramified $\mu$ and $\mu\nu^{-1}(\uf)\not =\abs{\uf}^{-1}$ if $\pi_{v}$ is unramified or special. In regards to Whittaker coefficients at ramified primes, we introduce the following.
 
 \begin{defn} 
 For $v \in B(\chi)$, let $\Psi_{v}:K_{v}^\x\to\C^\x$ be a character given by
\beq\label{E:char.W}\Psi_{v}(t):=\mu(\rmN_{v}(t))\cdot\chi\Abs_{K_{v}}^\onehalf(t)\eeq
for $N_{v}:K_{v}^\times \rightarrow \Q_{v}^\times$ the norm.
\end{defn}
\noindent Let $\mu_p(\Psi_{v})$ be a local invariant given by
\beq\label{E:9.W}\mu_p(\Psi_{v}):=\inf_{x\in K_v^\x}v_p(\Psi_{v}(x)-1).\eeq
Let $\cT\subset \A_K^\x$ be the subgroup consisting of ideles $z=(z_v)\in\prod_v K_v^\x$ with $z_v/\ol{z_v}\in\cO_{K_v}^\x$ for the primes $v$ split in $K$. 
Let $\rho: {\bf{A}}_{K}^\times \hookrightarrow \GL_{2}({\bf{A}})$ be a torus embedding and 
$\varsigma \in \GL_{2}({\bf{A}}_{K}^{(\ell p)})$ as in \S3. 
For $a \in ({\bf{A}}_{K}^{(\ell p)})^\times$, let
$
\rho_{\varsigma}(a)= \varsigma^{-1}\rho(a) \varsigma.
$
\\
\begin{defn}[Toric modular forms]\label{D:toric} We say that 
a nearly holomorphic form $f_{\chi} \in \pi$
is toric of character $\chi$ if
\[f_{\chi}|\rho_\varsigma(t)=\chi^{-1}(t)\cdot f_{\chi}\text{ for all }t\in\cT_v.\]
\end{defn}
\begin{remark} 
The notion is partly motivated by Waldspurger formula. Let $(\pi,\chi)$ be a pair of cuspidal automorphic representation of $\GL_2(\A)$ and a Hecke character over $K$ such that the corresponding Rankin--Selberg convolution $\pi \times \chi$ is self-dual with root number $1$. Then, the central L-value $L(\frac{1}{2},\pi \times \chi)$ is non-vanishing if and only if the period of a toric modular form against $\chi$ over $K$ is non-vanishing (\cf \cite{W}).
\end{remark}

\subsection{Fourier expansion} 
In this subsection, we consider the Fourier expansion of toric modular forms.\\
\\
Let the notation and hypothesis be as in \S A.1. To begin with, we have the following regarding toric modular forms.

\begin{prop}\label{torex} 
Let $p$ be an odd prime and $N$ a positive integer such that $p \ndivide N$. 
Let $K$ be an imaginary quadratic field satisfying (ord) and (Hg) with $\mathfrak{p}$ a prime above $p$ determined via an initial embedding $\iota_p: \overline{\Q} \hookrightarrow \C_p$. 
Let $f\in S_{k}(\Gamma_{0}(N),\epsilon)$ be a normalised newform with $T_{p}$-eigenvalue $\bfa_{p}(f)$ for $k=r+2$ with $r \geq 0$ and $\chi \in \Sg_{r,r}^{(2)}(c,\mathfrak{N},\epsilon)$ a Hecke character with infinity type $(r+1-j,1+j)$, 
for some $c$ prime to $N$ and $0 \leq j \leq r$. 
Let $N=N^{+}N^-$ for $N^-$ precisely divisible by inert or ramified primes in the extension $K/\Q$ and 
suppose that 
\begin{itemize}
\item[] $N^-$ is square-free.
\end{itemize}
Then, the following holds.

\begin{itemize}
\item[(i).] There exists a normalised toric modular form $\bff_{\chi}^{*} \in \pi$ satisfying the condition (MC) (\cf \S3). 
Let $\frakc$ be a prime-to-$p$ ideal of $\Q$ and let $\bfc\in (\A_{K,f}^{\setp})^\x$ such that $\il_\Q(\bfc)=\frakc$. Then the Fourier expansion of $\bff_{\chi}^{*}$ at the cusp $(\Z,\frakc)$ is given by
\[\bff_{\chi}^{*}|_{(\Z,\frakc)}(q)=\sum_{\beta}\bfa_\beta(\bff_{\chi}^{*},\frakc)q^\beta,\]
where
\[\bfa_\beta(\bff_{\chi}^{*},\frakc)=
\beta^{k/2+j} \cdot
\prod_{v \in \bdh, v \neq p}\bfa_{\chi,v}^{*}(\beta_{v}\bfc^{-1}_{v}) \cdot
\chi_{\overline{\mathfrak{p}}}(\beta^{-1})\bbI_{\Z_{(p)}^{\times}} (\beta)
\] 
for local Whittaker coefficients $\bfa_{\chi,v}^{*}:\Q_{v}^{\times} \rightarrow \C$.
Moreover, $$\bff_{\chi}^*\in \bfN_{k+2j}^{j}(\lsgN,\EucO)$$ for a non-negative integer $n$ and $\EucO$ specified below.
\item[(ii).] Let $\EucO$ be the finite extension of $\cO_{E_{f}}$ generated by $\stt{\localK{\chi,v}^*(1)}_{v\in B(\chi)}$ and the values of $p$-adic avatar $\wh\chi$. Then, the following holds. 
\begin{mylist}\item
The local Fourier coefficients $\localK{\chi,v}^*$ take value in $\EucO$ for each finite place $v\ndivide p$. 
 \item If either $v\not\in B(\chi)$ is unramified or $v\in A(\chi)$, then $\localK{\chi,v}^*(1)=1$.  
 In the unramified case, $$\localK{\chi,v}^*(\uf)=W_v^0(\DII{\uf}{1})$$  
 for $W_{v}^{0}$ the new Whittaker function such that $W_{v}^{0}(1)=1$.
 
 \item If $v\ndivides N$ is ramified with $c_v(\chi)=0$, then $\localK{\chi,v}^*(\uf^{-1})=1$,
 \item If $v\in B(\chi)$, then $\mu_p(\Psi_{v})=0$ (\cf Definition A.1) if and only if there exists $\eta_v\in \Q_{v}^\x$ such that
\[\localK{\chi,v}^*(\eta_v)\not\con 0\pmod{\frakm_{p}}.\]
\end{mylist}

\end{itemize}
\end{prop}
\noindent For part (i) (resp. part (ii)), we refer to \cite[Lem. 3.2 and Lem. 3.17]{Hs3} (resp. \cite[Prop. 3.15]{Hs3}).\\
\begin{remark} 
The hypothesis $N^-$ being square-free plays a key role in part (ii)(4).
\end{remark}
\noindent We have the following key non-triviality of the Fourier coefficients of the toric modular form.\\
\begin{prop}
Let $(f,\chi)$ be as above. Let  $\bff_{\chi}^*$ be the normalised toric modular form as in Proposition \ref{torex}. Suppose that 
the residual Galois representation \[\rho_f|_{G_K}\pmod{\frakm_p}\text{ is absolutely irreducible},\]
\\
Then, $\bff_{\chi}^*$ satisfies the non-triviality hypothesis (H) (\cf \S3) if and only if
\[\mu_{p}(\Psi):=\sum_{v|\frakc^-_\chi}\mu_p(\Psi_{v})=0\]
(\cf Definition A.1.).

\end{prop}
\begin{proof} We describe sketch of the argument in \cite[\S6]{Hs3}.\\
\\
As $\bff_{\chi}^*$ is toric, we have 
$$(\bff_{\chi}^*)^{\mathcal{R}}=|\Delta^{alg}|\cdot \bff_{\chi}^*$$
by definition (\cite[Lem. 6.1]{Hs2}). Recall that $|\Delta^{alg}|$ is a power of two.\\
\\
From the formula of $\bfa_\beta(\bff^*_{\chi},\frakc(a))$ in Proposition \ref{torex}, we note
\[\mu_p(\Psi)>0\text{ for some $v|\frakc_\chi^-$}\imply \bfa_\beta(\bff^*_{\chi},\frakc(a))\con 0\pmod{\frakm_{p}}\text{ for all }a\in \AFf^\x.\]
\\
Conversely, suppose that $\mu_p(\Psi_v)=0$ for all $v|\frakc^-_\chi$ and that the non-triviality hypothesis (H) does not hold. 
For $a\in\A_{K,f}^{(p N)}$,  
Proposition \ref{torex} readily implies
\begin{align*}&\bfa_\beta(\bff_{\chi}^*,\frakc(a))\con 0\pmod{\frakm_{p}}\text{ for all }
\beta\in\Q_{>0}\\
\iff&\localK{\chi}^\setp(\beta \bfc^{-1}\rmN(a^{-1}))\con 0\pmod{\frakm_{p}}\text{ for all }\beta\in
\Z_{\setp}^\x\end{align*}
where $N:\A_{K}^\times \rightarrow \A^\times$ the norm map. 
Here the global prime-to-$p$ Fourier coefficient $\localK{\chi}^\setp:(\AFf^\setp)^\x\to\C$ is given by
\beq\label{E:FC.W}\begin{aligned}\localK{\chi}^\setp(a):=&\prod_{v\in\bdh,v\ndivide p}\localK{\chi,v}^*(a_v).
\end{aligned}
\eeq
As a function on $(\AFf^\setp)^\x$, we thus have
\beq\label{E:8.W}\begin{aligned}\localK{\chi}^\setp(a)&\con 0\pmod{\frakm_{p}}\text{ for all }\\
&a\in \Z_{\setp}^\x \bfc^{-1} \det(U(N))\rmN((\A_{K,f}^{(p N)})^\x)=\Q^\x\bfc^{-1}\rmN((\A_{K,f}^\setp)^\x)\end{aligned}\eeq
with $U(N)$ the underlying level.\\
\\
From part (ii) of Proposition \ref{torex}, there exists $\eta=(\eta_v)\in\prod_{v|\frakc^-_\chi}\Q^\x_v$ such that $\localK{\chi,v}^*(\eta_v)\not \con 0\pmod{\frakm_{p}}$ for each $v|\frakc_\chi^-$. 
We now extend $\eta$ to be the idele in $\AFf^\x$ such that $\eta_v=1$ at $v\ndivides\frakc_\chi^-$. 
In view of \eqref{E:8.W} and the factorization formula of $\localK{\chi}^\setp$ (\cf \eqref{E:FC.W}), 
we thus have
\beq\label{E:14.W}\begin{aligned}\localK{\chi}^\setp(\eta\uf_v)&\con 0\pmod{\frakm_{p}}\iff W_v^0(\DII{\uf_v}{1})\con 0\pmod{\frakm_{p}}\text{ whenever }\\
&\uf_v\in [\eta^{-1}\bfc^{-1}]:=\Q^\x\eta^{-1}\bfc^{-1}\rmN((\A_{K,f}^\setp)^\x).\end{aligned}\eeq 
Here $v\ndivides pN$ and also $v \notin B(\chi)$.\\
\\
New local Whittaker coefficients are related to the underlying Galois representation $\rho_f$ as follows. From local-global compatibility, we have
\[\Tr\rho_{f}(Frob_{v})=\epsilon(\uf_v)^{-1}\abs{\uf_v}^{-k/2}W_v^0(\DII{\uf_v}{1})\text{ for all }v\ndivide pN.\]
Let $\rec_{K/\Q}:\A_{K}^\x\to \Gal(K/\Q)$ be the surjection induced via the reciprocity law. In view of  \eqref{E:14.W}, the above equality thus leads to 
\begin{align*}\Tr\rho_{f}(Frob_v)&\con 0\pmod{\frakm_{p}}\text{ whenever }\\
& Frob_v|_{K}=\rec_{K/\Q}(\uf_v)=\rec_{K/\Q}(\eta^{-1}\bfc^{-1}).\end{align*}
We thus note that $\rec_{K/\Q}(\eta^{-1}\bfc^{-1})$ is the complex conjugation $c$. 
In view of \cite[Lem. 6.2]{Hs3}, it thus follows that the Galois representation $\rho_{f} \mod{\frakm_{p}}$ is reducible. This contradiction finishes the proof.\\
\end{proof}
\begin{remark} 
The vanishing of $\mu_{p}(\Psi)$ is implied by the hypotheses (3) in Theorem A.\\
\end{remark}
\noindent\\

\end{document}